# Topology Optimization Using Polytopes

Arun L. Gain[a], Glaucio H. Paulino[a,1], Leonardo Duarte[b], Ivan F.M. Menezes[b]

[a] *Department of Civil & Environmental Engineering, University of Illinois at Urbana-Champaign, USA*

[b] *Tecgraf, Pontifical Catholic University of Rio de Janeiro (PUC-Rio), Brazil*

[1] *Corresponding author:* paulino@illinois.edu

---

**Abstract**

Meshing complex engineering domains is a challenging task. Arbitrary polyhedral meshes can provide the much needed flexibility in automated discretization of such domains. The geometric property of the polyhedral meshes such as the unstructured nature and the facial connectivity between elements makes them specially attractive for topology optimization applications. Numerical anomalies in designs such as the single node connections and checkerboard pattern, which are difficult to manufacture physically, are naturally alleviated with polyhedrons. Special interpolants such as Wachspress, mean value coordinates, maximum entropy shape functions are available to handle arbitrary shaped elements. But the finite elements approaches based on these shape functions face some challenges such as accurate and efficient computation of the shape functions and their derivatives for the numerical evaluation of the weak form integrals. In the current work, we solve the governing three-dimensional elasticity state equation using a Virtual Element Method (VEM) approach. The main characteristic difference between VEM and standard finite element methods (FEM) is that in VEM the canonical basis functions are not constructed explicitly. Rather the stiffness matrix is computed directly utilizing a projection map which extracts the linear component of the deformation. Such a construction guarantees the satisfaction of the patch test (used by engineers as an indicator of optimal convergence of numerical solutions under mesh refinement). Finally, the computations reduce to the evaluation of matrices which contain purely geometric surface facet quantities. The present work focuses on the first-order VEM in which the degrees of freedom associated with the vertices. Utilizing polyhedral elements for topology optimization, we show that the mesh bias in the member orientation is alleviated. The features of the current approach are demonstrated using numerical examples for compliance minimization and compliant mechanism problems.

---

## 1. Introduction

Numerical simulation of a typical engineering problem often begins with meshing a complicated domain, which is a challenging task. In general, unstructured meshes are preferred as they have shown to produce reliable and more accurate numerical solutions [16]. Thus, the polyhedral elements have become increasingly popular due to the flexibility they impart in automatic discretization of complicated design domains. Also, physical attributes such as facial connectivity between neighboring polyhedrons and unstructured nature of the mesh, make them specially attractive in topology optimization. Due to their characteristic geometry, numerical anomalies such as single node connections and checkerboarding are naturally alleviated. Orthogonally intersecting tension and compression members, distinctive features of an optimum structural layout [30, 23], are conveniently captured by polyhedral meshes. Thus, the goal of the current work is to develop a numerically efficient and accurate topology optimization approach for arbitrary polyhedral elements.

In the past, special interpolants have been proposed for arbitrary shaped elements such as Wachspress [45, 47], Sibson coordinates [34, 38], non-Sibson coordinates or Laplace shape functions [8, 24, 14, 38, 39], mean value coordinates [18, 19], metric coordinate method [28], maximum entropy shape function [36, 2, 25]. A detailed overview of the main developments in the field of conforming polygonal interpolants is provided in [37]. A summary of the polygonal/polyhedral interpolants is shown in Table. 1.

For two-dimensional topology optimization, polygonal shape functions have been explored in the past [42, 44, 21]. Utilizing the iso-parametric mapping scheme for numerically integration, polygonal shape functions can be efficiently implemented [39]. Polygonal shape functions and their derivatives are computed once for different reference n-gons and stored. Subsequently, these quantities can be retrieved as and when required. In three-dimensions, however, such a mapping scheme does not exist. Thus, numerical integration for each element are performed in physical coordinates which is computationally expensive. Moreover, to achieve accurate results, a high order numerical quadrature is required. This difficulty arises due to the non-polynomial nature of these shape functions. Figure 1 illustrates the complexity of elements in a typical polyhedral mesh.

Recently developed Virtual Element Method addresses some of the challenges facing the polyhedral elements for three-dimensional applications. The VEM originated from Mimetic Finite Difference (MFD) methods which have been successfully applied to diffusion [12], fluid flow [7] and elasticity problems [4]. The MFDs differ from the standard finite element approaches in the sense that in MFDs there are no explicitly defined shape functions associated with the discrete degrees of freedom. Thus, the continuous differential operators such as the div, grad, curl and trace, are approximated or *mimicked* by their discrete counterparts which utilize the discrete quantities defined only at the degrees of freedom. This provides greater flexibility in the geometric shapes of the admissible elements. In fact, high quality, skewed, degenerate and even non-convex polyhedra are all admissible. The MFD methods have been evolving from a finite difference/finite volume framework towards a Galerkin finite element-type framework called the Virtual Element Method [6, 22]. The feature which distinguishes VEM from classical finite element methods is that VEM does not require explicit computation of the approximation space. The construction of the discrete bilinear form for the elasticity problem begins with the kinematic decomposition of the element deformation space into constant strain and higher-order modes. The linear deformations are captured ex-actly, thereby passing the engineering patch test. A projection map is defined which enables such a decomposition. Using the discrete version of the projection map, the element stiffness matrix is constructed which requires computation of surface integral of the canonical basis functions over the element boundary faces. In this work, we explore the effectiveness of first-order VEM[1] [22] for three-dimensional topology optimization considering linear elastic behavior.

**Table 1:** Summary of the polygonal/polyhedral interpolants.

| Interpolants | References | Concave elements | Remarks |
|---|---|---|---|
| Wachspress | [45, 47, 46] | No | Earliest interpolant based on rational polynomials |
| Sibson | [34, 38] | No | Utilizes Voronoi tessellations to construct interpolant which reduces to the ratios of areas of Voronoi cells |
| Non-Sibson | [8, 24, 14, 38] | No | Also based on Voronoi tessellations. The interpolant is a function of Lebesgue measure of Voronoi edge and $L^2$ distance norm |
| Mean value coordinates | [18, 19] | Yes | Interpolant is a function of geometric quantities - $L^2$ distance norm and area |
| Metric coordinate method | [28] | Yes | Uses measures such as edge length, signed area of triangle, and trigonometric functions of sine and cosine to construct the shape functions |
| Maximum entropy | [36, 2, 25] | Yes | Shape functions and their derivatives are obtained by maximizing the Shannon's entropy function under prescribed boundary conditions |
| Harmonic coordinates | [26, 29, 11] | Yes | Shape functions and their derivatives are obtained by solving the Laplace equation hierarchically |

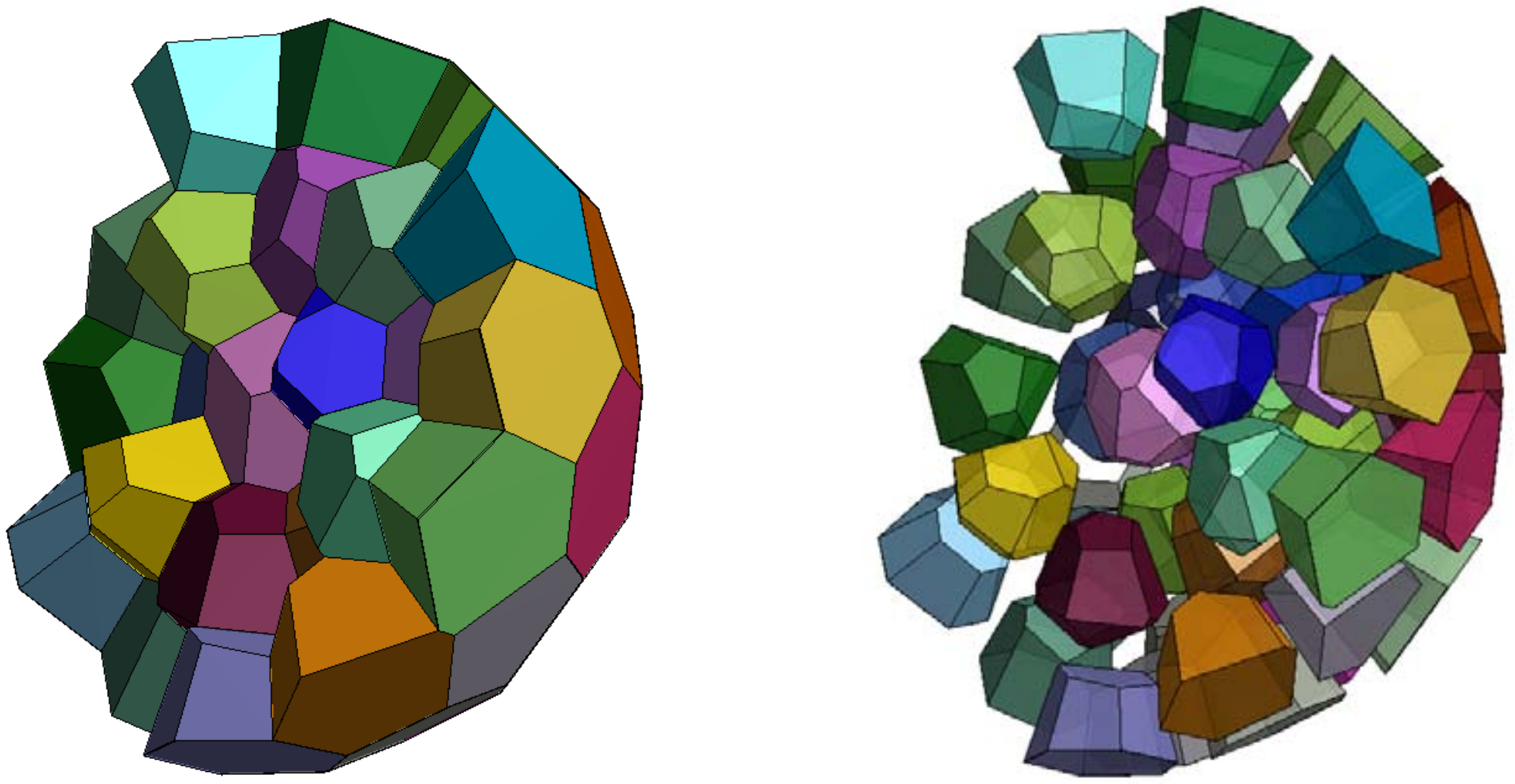

**Figure 1:** Illustration of the complexity of elements in an unstructured polyhedral mesh. (a) A section of the mesh for spherical domain. (b) Split view of the mesh.

The remainder of this paper is organized as follows. In Section 2, we define the governing elasticity problem and discuss the topology optimization problem formulation. Section 3 provides a brief overview of the Virtual Element Method for linear elasticity. A centroidal Voronoi tessellation (CVT) based meshing approach used to generated the polyhedral meshes is discussed in Section 4. In Section 5, we explore several numerical examples to evaluate our current approach. Finally, we provide some concluding remarks in Section 6.

## 2. Topology optimization problem formulation

In this work, we concentrate on topology optimization of linearized elastic system under small deformations subjected to surface tractions $\boldsymbol{t}$. The elasticity problem, for a smooth bounded domain $\Omega \subseteq \mathbb{R}^3$, is defined as follows. Find $\boldsymbol{u}$:

$$a(\boldsymbol{u}, \boldsymbol{v}) = f(\boldsymbol{v}), \qquad \forall \boldsymbol{u}, \boldsymbol{v} \in \boldsymbol{\mathcal{V}} \tag{1}$$

[1]Polyhedral elements are considered to be linear with three degrees of freedom associated with each vertex.

where

$$a(\boldsymbol{u},\boldsymbol{v}) = \int_\Omega \boldsymbol{C}\boldsymbol{\epsilon}(\boldsymbol{u}) : \boldsymbol{\epsilon}(\boldsymbol{v})\, d\boldsymbol{x}, \qquad f(\boldsymbol{v}) = \int_{\Gamma_t} \boldsymbol{t}\cdot\boldsymbol{v}\, d\boldsymbol{s}$$
$$\boldsymbol{\mathcal{V}} = \left\{\boldsymbol{v} \in H^1(\Omega)^3 : \boldsymbol{v}|_{\Gamma_u} = \boldsymbol{0}\right\} \tag{2}$$

Here, $\boldsymbol{\epsilon}(\boldsymbol{u}) = {}^1/_2(\nabla\boldsymbol{u} + \nabla\boldsymbol{u}^T)$ is the second-order linearized strain tensor and $\boldsymbol{C}$ is the elasticity tensor. In topology optimization, the working domain, $\Omega$, contains all the admissible shapes $\omega$, i.e., $\omega \subseteq \Omega$. Its boundary $\partial\Omega$ consists of three disjoint segments, $\partial\Omega = \Gamma_u \cup \Gamma_{t0} \cup \Gamma_t$, where $\Gamma_u$, $\Gamma_{t0}$, and $\Gamma_t$ represent displacement, homogeneous traction, and non-homogeneous traction boundary conditions ($\boldsymbol{t} \neq \boldsymbol{0}$), respectively. Also, the design $\omega$, with boundary $\partial\omega = \gamma_u \cup \gamma_{t0} \cup \gamma_t$, is constrained to satisfy $\gamma_u \subseteq \Gamma_u$ and $\gamma_t = \Gamma_t$. Here, $\gamma_u$, $\gamma_{t0}$, and $\gamma_t$ correspond to the boundaries of $\omega$ with displacement, homogeneous traction, and non-homogeneous traction boundary conditions, respectively (c.f. Fig 2). Body forces are ignored. The VEM is used to solve the state equation (1).

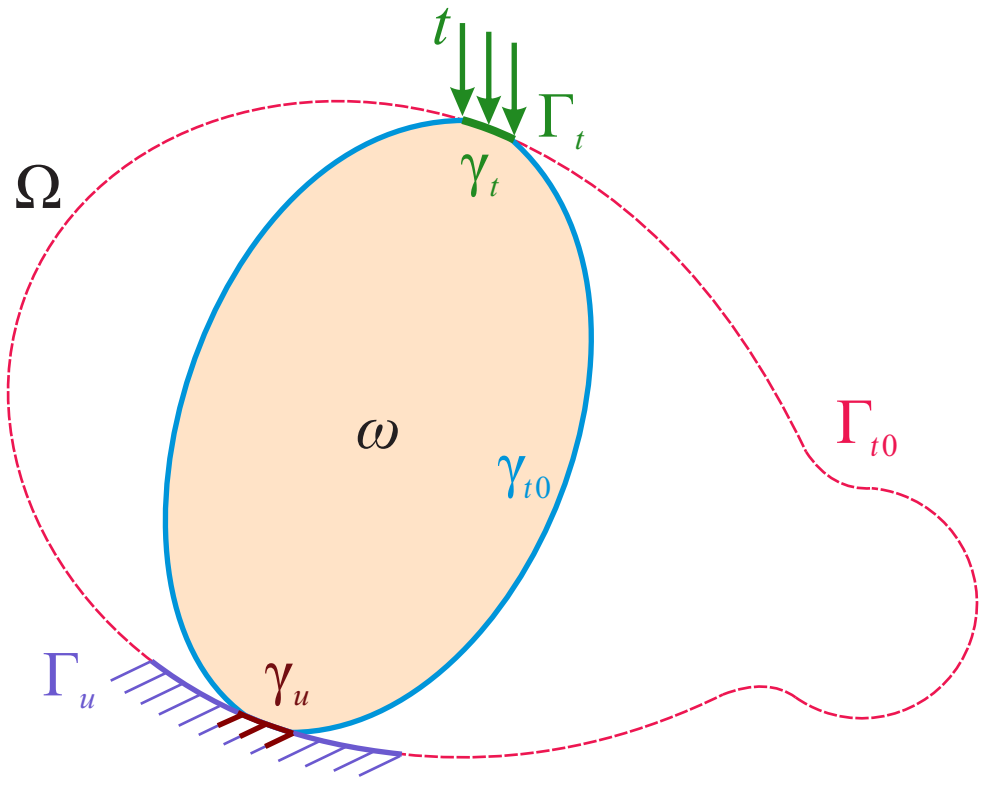

**Figure 2:** Domain description for the topology optimization problem. The boundary, $\partial\Omega$, of the working domain, $\Omega$, consists of $\Gamma_u$ (displacement boundary), $\Gamma_{t0}$ (homogeneous traction boundary) and $\Gamma_t$ (non-homogeneous traction boundary). The design $\omega$, with boundary $\partial\omega = \gamma_u \cup \gamma_{t0} \cup \gamma_t$, is constrained to satisfy $\gamma_u \subseteq \Gamma_u$ and $\gamma_t = \Gamma_t$. Boundaries $\gamma_u$, $\gamma_{t0}$, and $\gamma_t$ correspond to displacement, homogeneous traction, and non-homogeneous traction boundary conditions on $\partial\omega$, respectively.

We shall focus on two categories of topology optimization problems, compliance minimization and linear compliant mechanism. The topology optimization problem of compliance minimization refers to finding the stiffest configuration under applied loads and boundary conditions. Compliance, the work done by the loads, is defined as:

$$J_1(\rho) = \int_{\Gamma_t} \boldsymbol{t}\cdot\boldsymbol{u}\, d\boldsymbol{s} = \int_\Omega \boldsymbol{C}(\rho)\,\boldsymbol{\epsilon}(\boldsymbol{u}) : \boldsymbol{\epsilon}(\boldsymbol{u})\, d\boldsymbol{x} \tag{3}$$

The effective elasticity tensor $\boldsymbol{C}$ for the domain $\Omega$ is a function of density $\rho(\boldsymbol{x})$ and, as per the Solid Isotropic Material with Penalization (SIMP) model [9, 33], is expressed as:

$$\boldsymbol{C}(\rho) = [\varepsilon + (1-\varepsilon)\,\rho^p]\,\boldsymbol{C}^0 \tag{4}$$

The solid and void regions are filled with material of elasticity tensor $\boldsymbol{C}^0$ and $\varepsilon\boldsymbol{C}^0$, respectively, where $\varepsilon$ is chosen as $10^{-4}$. The penalization parameter, $p$, is set to 3.

The second category of problems we study are the linear compliant mechanisms, specifically the displacement inverter and the gripper problem. The objective is to maximize the displacement in a predefined direction, $u_{out}$, in response to the force, $f_{in}$, exerted by the actuator, modeled by a spring of stiffness $k_{in}$. So, the quantity we aim to minimize is:

$$J_2(\rho) = -u_{out} \tag{5}$$

Thus, the combined optimization problem, (3), (5), can be expressed as:

$$\begin{aligned} &\inf_{\rho} J_i(\rho) \qquad \text{for } i = 1, 2 \\ &\text{subject to:} \quad a(\boldsymbol{u}, \boldsymbol{v}) = f(\boldsymbol{v}), \qquad \int_{\Omega} \rho(\boldsymbol{x})\, d\boldsymbol{x} \leq V_f\, |\Omega| \end{aligned} \tag{6}$$

where $V_f$ is the prescribed maximum volume fraction and $|\cdot|$ denotes the measure (area or volume) of a set as well as the Euclidean norm of vector. In this work, we shall denote the components of vectors, matrices and tensors in the canonical Euclidean basis with subscripts inside parentheses (e.g. $\boldsymbol{u}_{(i)}$ or $\boldsymbol{\epsilon}_{(ij)}$) in order to make a distinction with indexed quantities.

For our simulations, we use the Optimality Criteria (OC) [10] as the optimization algorithm. Other mathematical programming algorithms such as, Method of Moving Asymptotes (MMA) [41], Sequential Linear Programming (SLP), Sequential Quadratic Programming (SQP), and CONvex LINearization approximations (CONLIN) [17] have also been used in density-based methods.

## 3. Virtual Element Method for linear elasticity

In this section, we present the Virtual Element Method for linear elasticity [22, 6] which is the governing state equation for the current topology optimization implementation. We discuss the mathematical framework of the method and also provide implementation details.

### 3.1. Theoretical background

Consider a partition of $\Omega$ into disjoint non-overlapping polyhedrons, $e$, of maximum diameter $h$. The Galerkin approximation $\boldsymbol{u}_h$ of $\boldsymbol{u}$, obtained by solving the discrete counterpart of (1), belongs to the conforming discrete space $\boldsymbol{\mathcal{V}}_h \subseteq \boldsymbol{\mathcal{V}}$ which consists of continuous displacement fields whose restriction to polyhedron $e$ belongs to the finite-dimensional space of smooth functions $\boldsymbol{\mathcal{W}}$. Space $\boldsymbol{\mathcal{W}}$ contains all the deformation states represented by the element $e$ - linear deformations and higher-order modes. Due to the conformity of $\boldsymbol{\mathcal{V}}_h$, the continuous $(a)$ and discrete bilinear form $(a_h)$ can be expressed as the corresponding sums of element contributions:

$$a(\boldsymbol{u}, \boldsymbol{v}) = \sum_e a^e(\boldsymbol{u}, \boldsymbol{v}), \qquad a_h(\boldsymbol{u}, \boldsymbol{v}) = \sum_e a_h^e(\boldsymbol{u}, \boldsymbol{v}) \tag{7}$$

In the space $\boldsymbol{\mathcal{W}}$, three degrees of freedom are associated with each vertex of a polyhedron. For this purpose, we consider the canonical basis $\boldsymbol{\varphi}_1, \ldots, \boldsymbol{\varphi}_{3n}$ of the form

$$\boldsymbol{\varphi}_{3i-2} = [\varphi_i,\ 0,\ 0]^T, \quad \boldsymbol{\varphi}_{3i-1} = [0,\ \varphi_i,\ 0]^T, \quad \boldsymbol{\varphi}_{3i} = [0,\ 0,\ \varphi_i]^T, \quad i = 1, \ldots, n \tag{8}$$

where $\varphi_1, \ldots, \varphi_n$ are a set of *barycentric coordinates*, for instance [45, 38, 18, 36, 2], which satisfy all the desired properties of a conforming interpolants such as partition of unity, Kronecker-delta, linear completeness and piece-wise linear ($C^0$ function) along the edges of $e$. We shall see later that VEM concerns only with the behavior of $\boldsymbol{\mathcal{W}}$ along the boundaries of the element $e$.

The VEM construction of the stiffness matrix begins with the kinematic decomposition of deformation states of $\boldsymbol{\mathcal{W}}$. To facilitate the derivations, the mean of the values of a function $\boldsymbol{v}$ sampled at the vertices of $e$ and volume average of $\boldsymbol{v}$ are represented by $\overline{\boldsymbol{v}}$ and $\langle \boldsymbol{v} \rangle$, respectively. We define the bases that span the space of linear deformations, $\boldsymbol{\mathcal{P}}$, over element $e$ as:

$$\begin{aligned}
\boldsymbol{p}_1(\boldsymbol{x}) &= \boldsymbol{e}_1, & \boldsymbol{p}_2(\boldsymbol{x}) &= \boldsymbol{e}_2, & \boldsymbol{p}_3(\boldsymbol{x}) &= \boldsymbol{e}_3, \\
\boldsymbol{p}_4(\boldsymbol{x}) &= (\boldsymbol{e}_1 \otimes \boldsymbol{e}_1)(\boldsymbol{x} - \overline{\boldsymbol{x}}), & \boldsymbol{p}_5(\boldsymbol{x}) &= (\boldsymbol{e}_1 \otimes \boldsymbol{e}_2)(\boldsymbol{x} - \overline{\boldsymbol{x}}), & \boldsymbol{p}_6(\boldsymbol{x}) &= (\boldsymbol{e}_1 \otimes \boldsymbol{e}_3)(\boldsymbol{x} - \overline{\boldsymbol{x}}), \\
\boldsymbol{p}_7(\boldsymbol{x}) &= (\boldsymbol{e}_2 \otimes \boldsymbol{e}_1)(\boldsymbol{x} - \overline{\boldsymbol{x}}), & \boldsymbol{p}_8(\boldsymbol{x}) &= (\boldsymbol{e}_2 \otimes \boldsymbol{e}_2)(\boldsymbol{x} - \overline{\boldsymbol{x}}), & \boldsymbol{p}_9(\boldsymbol{x}) &= (\boldsymbol{e}_2 \otimes \boldsymbol{e}_3)(\boldsymbol{x} - \overline{\boldsymbol{x}}), \\
\boldsymbol{p}_{10}(\boldsymbol{x}) &= (\boldsymbol{e}_3 \otimes \boldsymbol{e}_1)(\boldsymbol{x} - \overline{\boldsymbol{x}}), & \boldsymbol{p}_{11}(\boldsymbol{x}) &= (\boldsymbol{e}_3 \otimes \boldsymbol{e}_2)(\boldsymbol{x} - \overline{\boldsymbol{x}}), & \boldsymbol{p}_{12}(\boldsymbol{x}) &= (\boldsymbol{e}_3 \otimes \boldsymbol{e}_3)(\boldsymbol{x} - \overline{\boldsymbol{x}}).
\end{aligned} \tag{9}$$

Next, we define projection map $\pi_{\boldsymbol{\mathcal{P}}} : \boldsymbol{\mathcal{W}} \rightarrow \boldsymbol{\mathcal{P}}$ to extract linear deformations as:

$$\pi_{\boldsymbol{\mathcal{P}}} \boldsymbol{v} = \overline{\boldsymbol{v}} + \langle \nabla \boldsymbol{v} \rangle (\boldsymbol{x} - \overline{\boldsymbol{x}}) \tag{10}$$

We observe that the volumetric integral $\langle \nabla \boldsymbol{v} \rangle$ can be converted as:

$$\langle \nabla \boldsymbol{v} \rangle = \frac{1}{|e|} \int_e \nabla \boldsymbol{v}\ d\boldsymbol{x} = \frac{1}{|e|} \int_{\partial e} \boldsymbol{v} \otimes \boldsymbol{n}\ d\boldsymbol{s}, \tag{11}$$

An important property of $\pi_{\boldsymbol{\mathcal{P}}}$ which ensures consistency of the VEM approach is that $(\boldsymbol{v} - \pi_{\boldsymbol{\mathcal{P}}} \boldsymbol{v})$ is *energetically* orthogonal to $\boldsymbol{\mathcal{P}}$, $\forall \boldsymbol{v} \in \boldsymbol{\mathcal{W}}$, i.e.

$$a^e(\boldsymbol{p}, \boldsymbol{v} - \pi_{\boldsymbol{\mathcal{P}}} \boldsymbol{v}) = 0, \qquad \forall \boldsymbol{p} \in \boldsymbol{\mathcal{P}},\ \boldsymbol{v} \in \boldsymbol{\mathcal{W}} \tag{12}$$

To verify this identity, we can use the fact that the stress $\boldsymbol{\sigma}(\boldsymbol{p})$ is constant and $\langle \boldsymbol{\epsilon}(\boldsymbol{v}) \rangle = \boldsymbol{\epsilon}(\pi_{\boldsymbol{\mathcal{P}}} \boldsymbol{v})$, i.e. the volume average of the strain of a field $\boldsymbol{v}$ is the same as the strain of its linear deformation component, $\pi_{\boldsymbol{\mathcal{P}}} \boldsymbol{v}$. It will be useful to express the linear projection map $\pi_{\boldsymbol{\mathcal{P}}}$ in terms of the bases of $\boldsymbol{\mathcal{P}}$ as:

$$\begin{aligned}
\pi_{\boldsymbol{\mathcal{P}}} \boldsymbol{v} = {} & (\overline{\boldsymbol{v}})_{(1)} \boldsymbol{p}_1 + (\overline{\boldsymbol{v}})_{(2)} \boldsymbol{p}_2 + (\overline{\boldsymbol{v}})_{(3)} \boldsymbol{p}_3 + \langle \nabla \boldsymbol{v} \rangle_{(11)} \boldsymbol{p}_4 + \langle \nabla \boldsymbol{v} \rangle_{(12)} \boldsymbol{p}_5 + \langle \nabla \boldsymbol{v} \rangle_{(13)} \boldsymbol{p}_6 \\
& + \langle \nabla \boldsymbol{v} \rangle_{(21)} \boldsymbol{p}_7 + \langle \nabla \boldsymbol{v} \rangle_{(22)} \boldsymbol{p}_8 + \langle \nabla \boldsymbol{v} \rangle_{(23)} \boldsymbol{p}_9 + \langle \nabla \boldsymbol{v} \rangle_{(31)} \boldsymbol{p}_{10} + \langle \nabla \boldsymbol{v} \rangle_{(32)} \boldsymbol{p}_{11} + \langle \nabla \boldsymbol{v} \rangle_{(33)} \boldsymbol{p}_{12}
\end{aligned} \tag{13}$$

Finally, using $\pi_{\mathcal{P}}$, any deformation state $\boldsymbol{v} \in \boldsymbol{\mathcal{W}}$ can be kinematically decomposed as:

$$\boldsymbol{v} = \pi_{\mathcal{P}}\boldsymbol{v} + (\boldsymbol{v} - \pi_{\mathcal{P}}\boldsymbol{v}) \tag{14}$$

The higher-order component $(\boldsymbol{v} - \pi_{\mathcal{P}}\boldsymbol{v})$ belongs to a $(3n - 12)$ dimensional subspace of $\boldsymbol{\mathcal{W}}$.

Based on the kinematic decomposition of deformation state (14) and the energy orthogonality property (12), the continuous bilinear form can be written as:

$$a^e(\boldsymbol{u}, \boldsymbol{v}) = a^e(\pi_{\mathcal{P}}\boldsymbol{u}, \pi_{\mathcal{P}}\boldsymbol{v}) + a^e(\boldsymbol{u} - \pi_{\mathcal{P}}\boldsymbol{u}, \boldsymbol{v} - \pi_{\mathcal{P}}\boldsymbol{v}) \tag{15}$$

The first term corresponds to the constant strain modes and can be computed exactly by knowing the volume of the element $e$ (note that the arguments of the bilinear term are linear). Thus, ensuring that the engineering patch test is passed. The second term, corresponding to higher-order deformation modes, is difficult to compute. We replace this term by a crude estimate $s^e$, which can be conveniently computed, without affecting the energy associated with the constant strain modes. Thus, the discrete bilinear form is defined as:

$$a^e(\boldsymbol{u}, \boldsymbol{v}) \doteq a^e(\pi_{\mathcal{P}}\boldsymbol{u}, \pi_{\mathcal{P}}\boldsymbol{v}) + s^e(\boldsymbol{u} - \pi_{\mathcal{P}}\boldsymbol{u}, \boldsymbol{v} - \pi_{\mathcal{P}}\boldsymbol{v}) \tag{16}$$

Here, $s^e$ is any symmetric positive definite bilinear form and is chosen such that it has similar energy as the consistency term, $s^e(\cdot,\cdot) = a^e(\cdot,\cdot)$. We make a computationally inexpensive choice as [5]:

$$s^e(\boldsymbol{u}, \boldsymbol{v}) = \sum_{i=1}^{n} \alpha^e \boldsymbol{u}(\boldsymbol{x}_i) \cdot \boldsymbol{v}(\boldsymbol{x}_i) \tag{17}$$

where $\alpha^e$ is a positive coefficient that ensures correct scaling of the energies of higher-order mode. In the following section we provide explicit expressions for the element stiffness matrix.

### *3.2. Implementation details*

Here we show the implementation details of the VEM formulation discussed in the previous section. We provide the expression for discrete linear projection map which is later used to construct the stiffness matrix terms. Let, $\boldsymbol{P}_{\mathcal{P}}$ be the discrete representation of the projection $\pi_{\mathcal{P}}$, i.e.,

$$\pi_{\mathcal{P}}\boldsymbol{\varphi}_j = \sum_{k=1}^{3n} (\boldsymbol{P}_{\mathcal{P}})_{(kj)}\, \boldsymbol{\varphi}_k \tag{18}$$

To obtain $\boldsymbol{P}_{\mathcal{P}}$, we use (13) to express $\pi_{\mathcal{P}}$ in terms of its corresponding bases as:

$$\pi_{\mathcal{P}}\boldsymbol{\varphi}_j = \sum_{\ell=1}^{12} (\boldsymbol{W}_{\mathcal{P}})_{(j\ell)}\, \boldsymbol{p}_\ell \tag{19}$$

where $\boldsymbol{W}_{\boldsymbol{\mathcal{P}}}$ is an $3n \times 12$ matrix whose $j$th row is given by,

$$
\begin{aligned}
&[\left(\overline{\boldsymbol{\varphi}_j}\right)_{(1)},\ \left(\overline{\boldsymbol{\varphi}_j}\right)_{(2)},\ \left(\overline{\boldsymbol{\varphi}_j}\right)_{(3)},\ \langle \nabla \boldsymbol{\varphi}_j \rangle_{(11)},\ \langle \nabla \boldsymbol{\varphi}_j \rangle_{(12)},\ \langle \nabla \boldsymbol{\varphi}_j \rangle_{(13)},\\
&\qquad \langle \nabla \boldsymbol{\varphi}_j \rangle_{(21)},\ \langle \nabla \boldsymbol{\varphi}_j \rangle_{(22)},\ \langle \nabla \boldsymbol{\varphi}_j \rangle_{(23)},\ \langle \nabla \boldsymbol{\varphi}_j \rangle_{(31)},\ \langle \nabla \boldsymbol{\varphi}_j \rangle_{(32)},\ \langle \nabla \boldsymbol{\varphi}_j \rangle_{(33)}]
\end{aligned} \tag{20}
$$

Using the linear precision property of the canonical basis functions (8), we express $\boldsymbol{p}_\ell$ in (19) in terms of its discrete counterpart (bases sampled at the vertices) and upon simplification, we obtain:

$$
\pi_{\boldsymbol{\mathcal{P}}} \boldsymbol{\varphi}_j = \sum_{k=1}^{3n} \left( \boldsymbol{N}_{\boldsymbol{\mathcal{P}}} \boldsymbol{W}_{\boldsymbol{\mathcal{P}}}^T \right)_{(kj)} \boldsymbol{\varphi}_k \tag{21}
$$

Subsequent comparison with (18), we get $\boldsymbol{P}_{\boldsymbol{\mathcal{P}}} = \boldsymbol{N}_{\boldsymbol{\mathcal{P}}} \boldsymbol{W}_{\boldsymbol{\mathcal{P}}}^T$, where $\boldsymbol{N}_{\boldsymbol{\mathcal{P}}}, \boldsymbol{W}_{\boldsymbol{\mathcal{P}}} \in \mathbb{R}^{3n \times 12}$. A block of three rows of $\boldsymbol{N}_{\boldsymbol{\mathcal{P}}}$, corresponding to $i$th vertex, is explicitly expressed as:

$$
\begin{aligned}
(\boldsymbol{N}_{\boldsymbol{\mathcal{P}}})_{(3i-2:3i,:)} = \Bigg[ &\begin{matrix} 1 & 0 & 0 & (\boldsymbol{x}_i - \overline{\boldsymbol{x}})_{(1)} & (\boldsymbol{x}_i - \overline{\boldsymbol{x}})_{(2)} & (\boldsymbol{x}_i - \overline{\boldsymbol{x}})_{(3)} \\ 0 & 1 & 0 & 0 & 0 & 0 \\ 0 & 0 & 1 & 0 & 0 & 0 \end{matrix}\\
&\begin{matrix} 0 & 0 & 0 & 0 & 0 & 0 \\ (\boldsymbol{x}_i - \overline{\boldsymbol{x}})_{(1)} & (\boldsymbol{x}_i - \overline{\boldsymbol{x}})_{(2)} & (\boldsymbol{x}_i - \overline{\boldsymbol{x}})_{(3)} & 0 & 0 & 0 \\ 0 & 0 & 0 & (\boldsymbol{x}_i - \overline{\boldsymbol{x}})_{(1)} & (\boldsymbol{x}_i - \overline{\boldsymbol{x}})_{(2)} & (\boldsymbol{x}_i - \overline{\boldsymbol{x}})_{(3)} \end{matrix} \Bigg]
\end{aligned} \tag{22}
$$

and a block of three rows of $\boldsymbol{W}_{\boldsymbol{\mathcal{P}}}$ is simplified as:

$$
\begin{aligned}
(\boldsymbol{W}_{\boldsymbol{\mathcal{P}}})_{(3i-2:3i,:)} = \Bigg[ &\begin{matrix} 1/n & 0 & 0 & (\boldsymbol{q}_i)_{(1)} & (\boldsymbol{q}_i)_{(2)} & (\boldsymbol{q}_i)_{(3)} \\ 0 & 1/n & 0 & 0 & 0 & 0 \\ 0 & 0 & 1/n & 0 & 0 & 0 \end{matrix}\\
&\begin{matrix} 0 & 0 & 0 & 0 & 0 & 0 \\ (\boldsymbol{q}_i)_{(1)} & (\boldsymbol{q}_i)_{(2)} & (\boldsymbol{q}_i)_{(3)} & 0 & 0 & 0 \\ 0 & 0 & 0 & (\boldsymbol{q}_i)_{(1)} & (\boldsymbol{q}_i)_{(2)} & (\boldsymbol{q}_i)_{(3)} \end{matrix} \Bigg]
\end{aligned} \tag{23}
$$

where the surface integral vector $\boldsymbol{q}_i$ is given by:

$$
\boldsymbol{q}_i = \frac{1}{|e|} \int_{\partial e} \varphi_i \boldsymbol{n} \mathrm{d}\boldsymbol{s} \tag{24}
$$

To compute the surface integral, such as the one encountered in (24), we use first-order accurate nodal quadrature scheme which does not require the computation of canonical basis functions inside the elements. The approximation of the surface integration of a generic function $g$ over a face $F$ is defined as:

$$
\int_F g \ d\boldsymbol{s} \doteq \sum_{j \in F} w_j^F g(\boldsymbol{x}_j) \tag{25}
$$

where $\{\boldsymbol{x}_j : j \in F\}$ represent the vertices on the face $F$. The nodal weight $w_j^F$ for vertex $\boldsymbol{x}_j$ is taken as the area of the quadrilateral formed by $\boldsymbol{x}_j$, the centroid of the face $F$, and mid-points of the edges incident on $\boldsymbol{x}_j$ (Fig. 3).

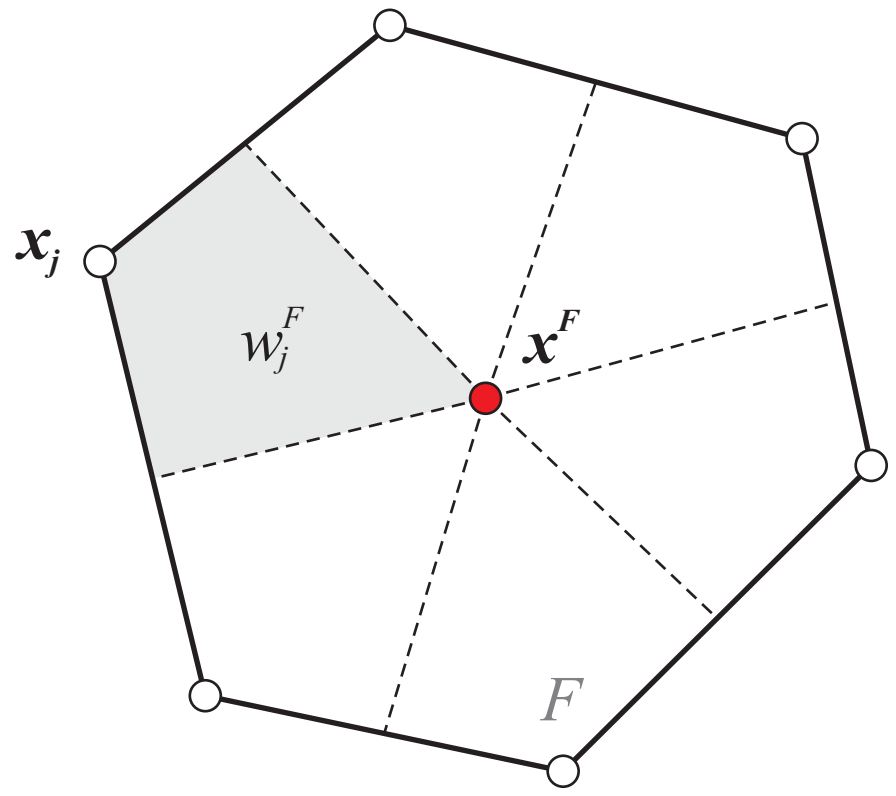

**Figure 3:** Nodal quadrature scheme for surface integration. The variables $\boldsymbol{x}^F$ and $w_j^F$ represent the centroid of the polygon and the nodal weight associated with the vertex $\boldsymbol{x}_j$ on any face $F$, respectively.

The element stiffness matrix $\boldsymbol{K}_e$ can now be obtained using the surface integral matrix $\boldsymbol{W}_{\mathcal{P}}$ and discrete projection $\boldsymbol{P}_{\mathcal{P}}$. From (16) we have,

$$(\boldsymbol{K}_e)_{(jk)} = a_h^e(\boldsymbol{\varphi}_j, \boldsymbol{\varphi}_k) = a^e(\pi_{\mathcal{P}}\boldsymbol{\varphi}_j, \pi_{\mathcal{P}}\boldsymbol{\varphi}_k) + s^e(\boldsymbol{\varphi}_j - \pi_{\mathcal{P}}\boldsymbol{\varphi}_j, \boldsymbol{\varphi}_k - \pi_{\mathcal{P}}\boldsymbol{\varphi}_k) \tag{26}$$

Using (19), the first term of $\boldsymbol{K}_e$ can be simplified as:

$$a^e(\pi_{\mathcal{P}}\boldsymbol{\varphi}_j, \pi_{\mathcal{P}}\boldsymbol{\varphi}_k) = |e| \left(\boldsymbol{W}_{\mathcal{P}}\boldsymbol{D}\boldsymbol{W}_{\mathcal{P}}^T\right)_{(jk)} \tag{27}$$

where the matrix $\boldsymbol{D}$ is a function of elasticity tensor $\boldsymbol{C}$ and is defined as:

$$(\boldsymbol{D})_{(\ell m)} = \frac{1}{|e|} a^e(\boldsymbol{p}_\ell, \boldsymbol{p}_m) = \boldsymbol{C}\boldsymbol{\epsilon}(\boldsymbol{p}_\ell) : \boldsymbol{\epsilon}(\boldsymbol{p}_m), \qquad \ell, m = 1, \ldots, 12 \tag{28}$$

Now, adopting the inexpensive choice of $s^e$ (17), the second term in the stiffness matrix is written as:

$$s^e(\boldsymbol{\varphi}_j - \pi_{\mathcal{P}}\boldsymbol{\varphi}_j, \boldsymbol{\varphi}_k - \pi_{\mathcal{P}}\boldsymbol{\varphi}_k) = \left[(\boldsymbol{I} - \boldsymbol{P}_{\mathcal{P}})^T \boldsymbol{S}^e (\boldsymbol{I} - \boldsymbol{P}_{\mathcal{P}})\right]_{(jk)} \tag{29}$$

where $(\boldsymbol{S}^e)_{(jk)} = s^e(\boldsymbol{\varphi}_j, \boldsymbol{\varphi}_k)$ and corresponds to $\boldsymbol{S}^e = \alpha^e \boldsymbol{I}_{3n}$. We need to ensure that $s^e(\cdot,\cdot)$ is of the same order of magnitude as $a^e(\cdot,\cdot)$, so an appropriate value of $\alpha^e$ is

$$\alpha^e = \bar{\alpha}^e \text{trace}(|e|\boldsymbol{W}_{\mathcal{P}}\boldsymbol{D}\boldsymbol{W}_{\mathcal{P}}^T) \tag{30}$$

Based on the studies conducted in reference [20], $\bar{\alpha}^e$ is chosen as 0.05. Finally, the element stiffness matrix is given by

$$\boldsymbol{K}_e = |e|\boldsymbol{W}_{\mathcal{P}}\boldsymbol{D}\boldsymbol{W}_{\mathcal{P}}^T + \alpha^e (\boldsymbol{I} - \boldsymbol{P}_{\mathcal{P}})^T (\boldsymbol{I} - \boldsymbol{P}_{\mathcal{P}}) \tag{31}$$

To compute the force vectors corresponding to a body force, a first-order accurate scheme, similar to the one discussed in (25), can be used. Nodal quadrature for body forces involves computing nodal volume weights for volume integration. The weight for vertex $\boldsymbol{x}_i$ is the volumes of the polyhedron formed by $\boldsymbol{x}_i$, the centroid of element $e$, the centroid of faces and the centroids of the edges incident on $\boldsymbol{x}_i$. In case of surface traction, (25) can be used. The global stiffness matrix and global force vector are obtain by standard assembly process. For a detailed study on the accuracy and convergence of VEM for three-dimensional elasticity refer to [22].

## 4. On centroidal Voronoi tessellation meshing

We use Voronoi diagram to generate three-dimensional polyhedral meshes [3, 48, 43]. Given a set of $n_s$ distinct seeds $\boldsymbol{\mathcal{S}} = \{\boldsymbol{s}_i\}_{i=1}^{n_s}$, the Voronoi tessellation of the domain $\Omega \subseteq \mathbb{R}^3$ is defined as:

$$\mathfrak{D}(\boldsymbol{\mathcal{S}}) = \{e \cap \Omega : \boldsymbol{s}_i \in \boldsymbol{\mathcal{S}}\} \tag{32}$$

where $e$ is the Voronoi cell corresponding to seed $\boldsymbol{s}_i$:

$$e = \left\{\boldsymbol{x} \in \mathbb{R}^3 : |\boldsymbol{x} - \boldsymbol{s}_i| < |\boldsymbol{x} - \boldsymbol{s}_j|, \quad \forall j \neq i\right\} \tag{33}$$

The above definition of $e$ represents a domain consisting of all points that are closer to seed $\boldsymbol{s}_i$ than any other seed $\boldsymbol{s}_j \in \boldsymbol{\mathcal{S}}$. Note that the Voronoi cells are necessarily convex polyhedrons since they are formed by the finite intersection of convex half-planes.

Following the guidelines discussed in [48, 43], a polyhedral discretization is obtained from the Voronoi diagram of a given set of seeds and their reflections about the closest boundary of $\Omega$, i.e. $\Gamma$. Our meshing algorithm is implemented for general domains using the concept of a signed distance function. A signed distance function $d_\Omega(\boldsymbol{x})$ is defined as:

$$d_\Omega(\boldsymbol{x}) = z(\boldsymbol{x}) \min\left(|\boldsymbol{x} - \boldsymbol{y}|\right), \qquad \forall \boldsymbol{y} \in \Gamma \tag{34}$$

where $z(\boldsymbol{x})$ is the sign function defined as:

$$z(\boldsymbol{x}) = \begin{cases} -1, & \boldsymbol{x} \in \Omega, \\ +1, & \boldsymbol{x} \in \mathbb{R}^3 \backslash \Omega. \end{cases} \tag{35}$$

Thus, $d_\Omega(\boldsymbol{x}) = 0$ if $\boldsymbol{x} \in \Gamma$ and $d_\Omega(\boldsymbol{x}) < 0$ if $\boldsymbol{x} \in \Omega \backslash \Gamma$. Using the signed distance function and its gradient, the reflection, $\boldsymbol{s}_i^r$, of the seed $\boldsymbol{s}_i$ can be calculated as:

$$\boldsymbol{s}_i^r = \boldsymbol{s}_i^r - 2 d_\Omega(\boldsymbol{s}_i) \nabla d_\Omega(\boldsymbol{s}_i) \tag{36}$$

First, to construct a polyhedron discretization of the domain $\Omega$, each point in $\boldsymbol{\mathcal{S}}$ is reflected about the closest boundary $\Gamma$. The resulting set of points are denoted by $\boldsymbol{\mathcal{R}}_\Omega(\boldsymbol{\mathcal{S}})$. Subsequently, we construct the Voronoi diagram of the space using the original point set and its reflection, $\boldsymbol{\mathcal{T}}(\boldsymbol{\mathcal{S}} \cup \boldsymbol{\mathcal{R}}_\Omega(\boldsymbol{\mathcal{S}}); \mathbb{R}^3)$. Thus, the discretization of $\Omega$ is the collection of the Voronoi cells associated with seeds $\boldsymbol{\mathcal{S}}$. For a given point set $\boldsymbol{\mathcal{S}}$, the discretization of the domain $\Omega$ is uniquely defined and denoted by:

$$\boldsymbol{\mathcal{M}}_\Omega(\boldsymbol{\mathcal{S}}) = \left\{ e \in \boldsymbol{\mathcal{T}}(\boldsymbol{\mathcal{S}} \cup \boldsymbol{\mathcal{R}}_\Omega(\boldsymbol{\mathcal{S}}); \mathbb{R}^3) : \boldsymbol{s}_i \in \boldsymbol{\mathcal{S}} \right\} \tag{37}$$

If the Voronoi cell of a seed $\boldsymbol{s}_i$ and its reflection have a common edge, then this edge forms an approximation to the domain boundary and a reasonable discretization of $\Omega$ is obtained. In order to mesh complicated geometries, a signed distance function along with set operations such as union, difference, and intersection are used [32, 43].

In our meshing algorithm, a set of signed distance functions corresponding to basic geometric shapes such as three-dimensional plane, sphere, cylinder and rectangular box is defined. We also construct a bounding box $B$, that contains the domain $\Omega$, to generate the random seeds in $\mathbb{R}^3$. A random seed is accepted only if it lies inside the domain $\Omega$, determined by evaluating the sign of the resultant distance function, $d_\Omega$, associated with $\Omega$. Algorithm 1 shows the basic steps for obtaining a random point set of size $n_s$.

**Algorithm 1** Initial random seed placement

**input:** $B, n_s$ %% $B \supset \Omega \in \mathbb{R}^3$ *and* $n_s$ *is the number of seeds*
  set $\boldsymbol{\mathcal{S}} = \emptyset$
  **while** $|\boldsymbol{\mathcal{S}}| < n_s$ **do**
    generate random point $\mathbf{y} \in B$
    **if** $d_\Omega(\mathbf{y}) < 0$ **then**
      $\boldsymbol{\mathcal{S}} \leftarrow \boldsymbol{\mathcal{S}} \cup \{\mathbf{y}\}$
    **end if**
  **end while**
**output:** $\boldsymbol{\mathcal{S}}$

We handle *convex* and *non-convex* features of $\Omega$ by carefully choosing a set of seeds to be reflected w.r.t. the boundary. Reflection of a seed far from the boundary may land inside the domain or interfere with the reflection of another seed. Since the reflection of most of the seeds in the interior of the domain has no effect on the approximation of the boundary, we reflect only the seeds that are in a band near the boundary. A seed $\boldsymbol{s}_i \in \boldsymbol{\mathcal{S}}$ is reflected about boundary segment $\Gamma$ provide that:

$$|d_\Omega(\boldsymbol{s}_i)| < c \left( \frac{|\Omega|}{n_s} \right)^{1/3} \tag{38}$$

where $c$ is the proportionality constant, chosen to be greater than 1 to make the band size near the boundary larger than the average element volume.

Clearly, Voronoi meshes generated from random/quasi-random seeds may cause inconsistencies at the boundaries resulting in a poor approximation of the boundaries of the design domain. To introduce some regularity in the Voronoi meshes, we construct Centroidal Voronoi tessellations (CVT) using a modified Lloyd's algorithm [27]. For a large number of iterations, CVT cells tend to be uniform in size [15]. To generate CVT meshes, we replace seeds $\boldsymbol{\mathcal{S}}$ with centroids $\boldsymbol{\mathcal{S}}_c$ of the Voronoi cells. We compute the polyhedron centroid by partitioning it into tetrahedrons and evaluating the weighted mean of the centroids of the resulting tetrahedrons. The weights are the volumes of the tetrahedrons.

To construct the element stiffness matrix using the approach discussed in Section 3, along with vertices locations and element connectivity information, we need to know the vertices incident on each face of polyhedrons. The pseudo-code listed in Algorithm 2 summarizes our approach to obtain facial information in MATLAB. The inputs to the algorithm are the element connectivity and the coordinates of the vertices and seeds. For each seed of the mesh, the algorithm computes the convex hull of the set $\boldsymbol{V}$ (vertices in an element). The convex hull of $\boldsymbol{V}$ is a matrix $\boldsymbol{H}$ with as many rows as the number of triangular subdivision of faces present in the convex hull, and three columns containing the indexes of the vertices of the corresponding triangular subdivisions. Then, we iterate over all triangles and unite those that are co-planar to obtain the polygonal faces of the polyhedron. The resulting array, called $\boldsymbol{Elm}$, contains the faces and vertices of each polyhedral element. Figure 4 shows some of the sample meshes and their statistical information, obtained from the current algorithm.

**Algorithm 2** Construction of final mesh consisting of vertices, elements and faces in MATLAB.

**input:** $e, \boldsymbol{N}, \boldsymbol{\mathcal{S}}, n_s$ %% *Voronoi cells* $e$, *vertices coords* $N$, *seeds coords* $\mathcal{S}$
  $\boldsymbol{Elm} \leftarrow 0$ %% *initialize an array with size of* $n_s$
  **for** $i = 1$ to $n_s$ **do**
    let $\boldsymbol{V} = \boldsymbol{N}(e(i))$
    construct convex hull $\boldsymbol{H} \leftarrow \boldsymbol{\mathcal{H}}(\boldsymbol{V}; \mathbb{R}^3)$
    $m \leftarrow |\boldsymbol{H}|$ %% *number of triangular subdivisions* $m$ *obtained from convex hull*
    *ElementFaces* $\leftarrow 0$
    **for** $j = 1$ to $m$ **do**
      let $\boldsymbol{T} = e(\boldsymbol{H}(j))$ %% *vertices of a triangle*
      $\boldsymbol{T} = OrderVertices(\boldsymbol{N}, \boldsymbol{T}, \boldsymbol{\mathcal{S}}(i))$
      *ElementFaces* $\leftarrow \boldsymbol{T}$ %% *create/unite faces*
    **end for**
    $\boldsymbol{Elm}(i) \leftarrow$ *ElementFaces*
  **end for**
**output:** $\boldsymbol{Elm}$

## 5. Numerical examples

In this section, numerical examples are shown to demonstrate the effectiveness of the current approach. We start with the benchmark cantilever beam problem solved on a box domain using different mesh discretizations, followed by problems on non-Cartesian design domains. For our simulations, the Young's modulus, $E$, and Poisson's ratio, $\nu$, are taken as $10,000$ and 0.3, respectively. As mentioned earlier, Optimality Criteria (OC) is the optimizer of choice. The optimization is terminated when either the maximum of the change in element densities is less than 0.01 or the number of iterations exceed 300.

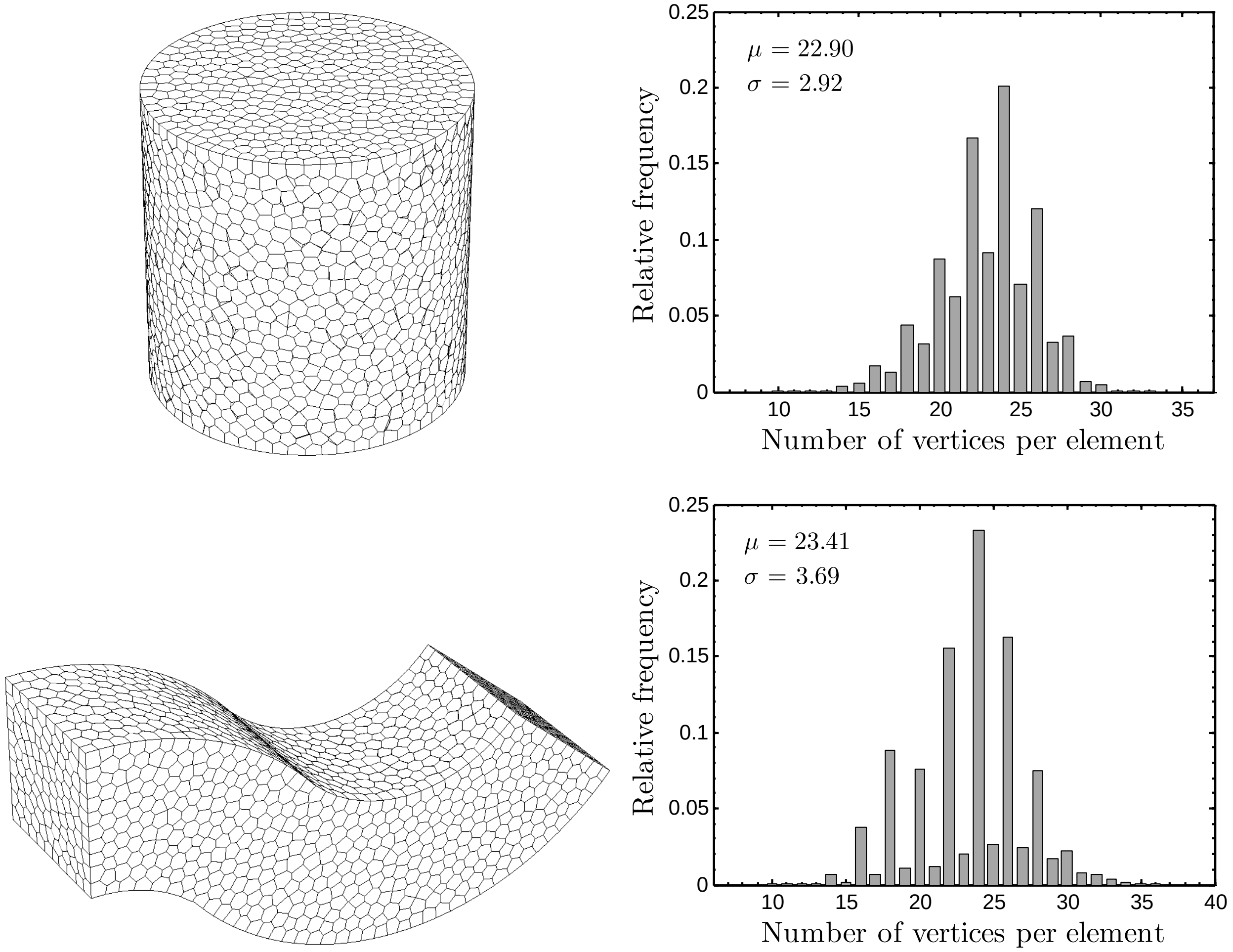

**Figure 4:** Sample meshes. (a) Cylinder. (b) Cylinder mesh statistics. On average, elements in the mesh have approximately 23 vertices with a standard deviation of 2.92. (c) Curved cantilever beam. (d) Curved cantilever mesh statistics. On average, the elements in the mesh also have approximately 23 vertices with a standard deviation of 3.69.

### *5.1. Cantilever beam problem on box domain*

We begin with the benchmark cantilever beam problem for a design domain of dimension $2 \times 1 \times 1$. The left face of the box is fixed and a point load is applied in the middle of the right face (refer Fig. 5). The problem is solved on both hexahedral and polyhedral element meshes. Taking advantage of symmetry, only half of the domain is optimized, which is discretized using 54,872 hexahedrons (60,060 nodes) and 10,000 polyhedrons (58,601 nodes). A linear filter of radius 5% of the maximum domain dimension is used and a volume fraction of 0.1 is prescribed.

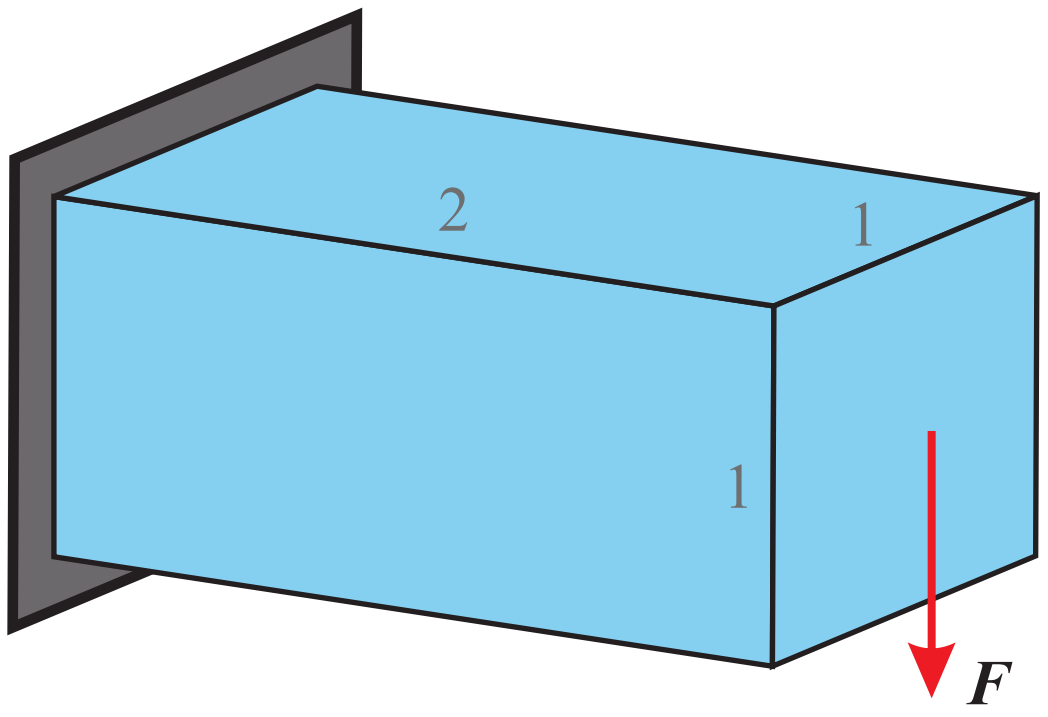

**Figure 5:** Cantilever beam problem.

With the present approach, both mesh discretizations produce similar optimization results (Fig. 6). Note that the optimized topologies shown in Fig. 6 (also all subsequent results) show only the elements whose density exceeds 0.5. The final compliance values are 0.1098 and 0.1082 for the topologies on hexahedral and polyhedral meshes, respectively.

For comparison, the same problem is solved using the finite element method on a hexahedral mesh. Topology similar to the present method is obtained. The convergence history for all three cases are illustrated in Fig. 7. As expected, a smooth monotonic convergence is obtained for all three cases and they all converge to similar final compliance values.

### *5.2. Shear loaded thin disc*

We next investigate the shear loaded thin disc problem. The thin disc domain has an external radius of 6 units with an internal cylindrical hole of radius 1 unit and has a thickness of 0.5 units (Fig. 8). Eight equidistant shear loads are applied along the circumference of the disc and all the nodes along the cylindrical hole are fixed. A polyhedral mesh of 10,000 elements (55,810 nodes) is used to discretize the design domain. A filter radius of 4% of the outer diameter is selected and a volume fraction of 0.2 is enforced.

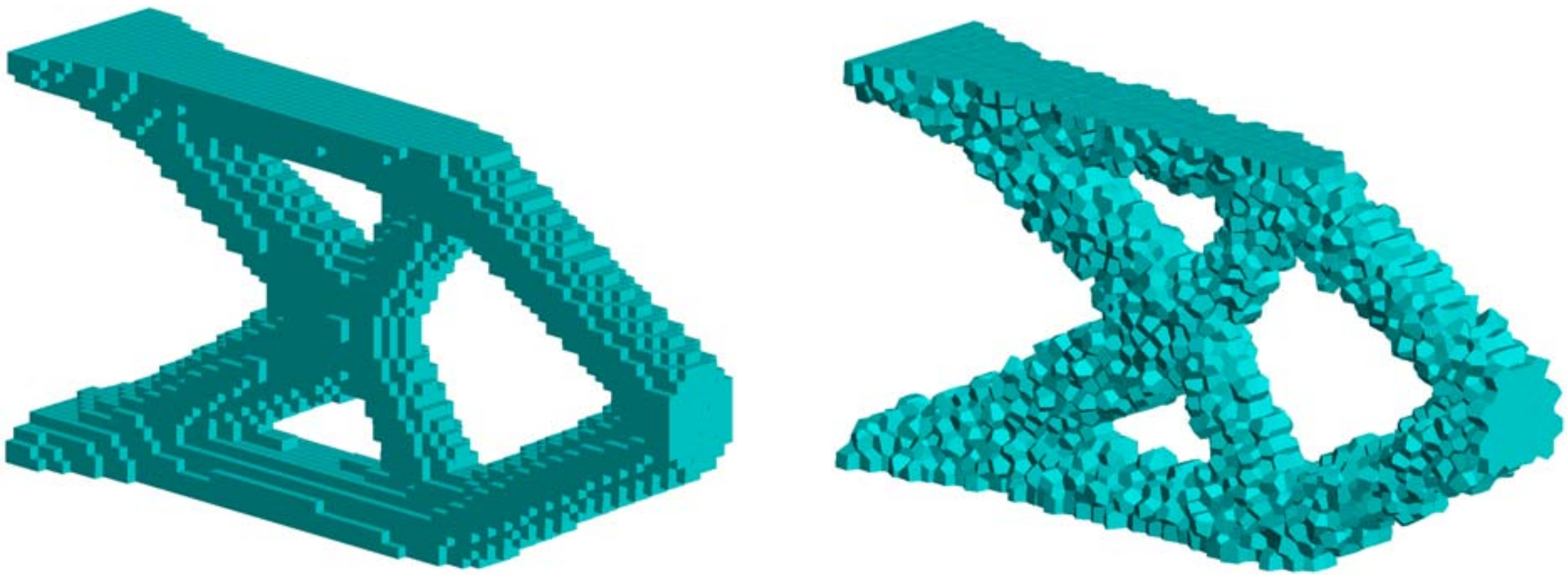

**Figure 6:** Converged topologies for the cantilever beam problem using the present method. (a) Hexahedral mesh of 54,872 elements, 60,060 nodes ($J = 0.1098$). (b) Polyhedral mesh of 10,000 elements, 58,601 nodes. The average number of vertices per polyhedron is, $\mu = 22.85$, with standard deviation, $\sigma = 3.80$, ($J = 0.1082$).

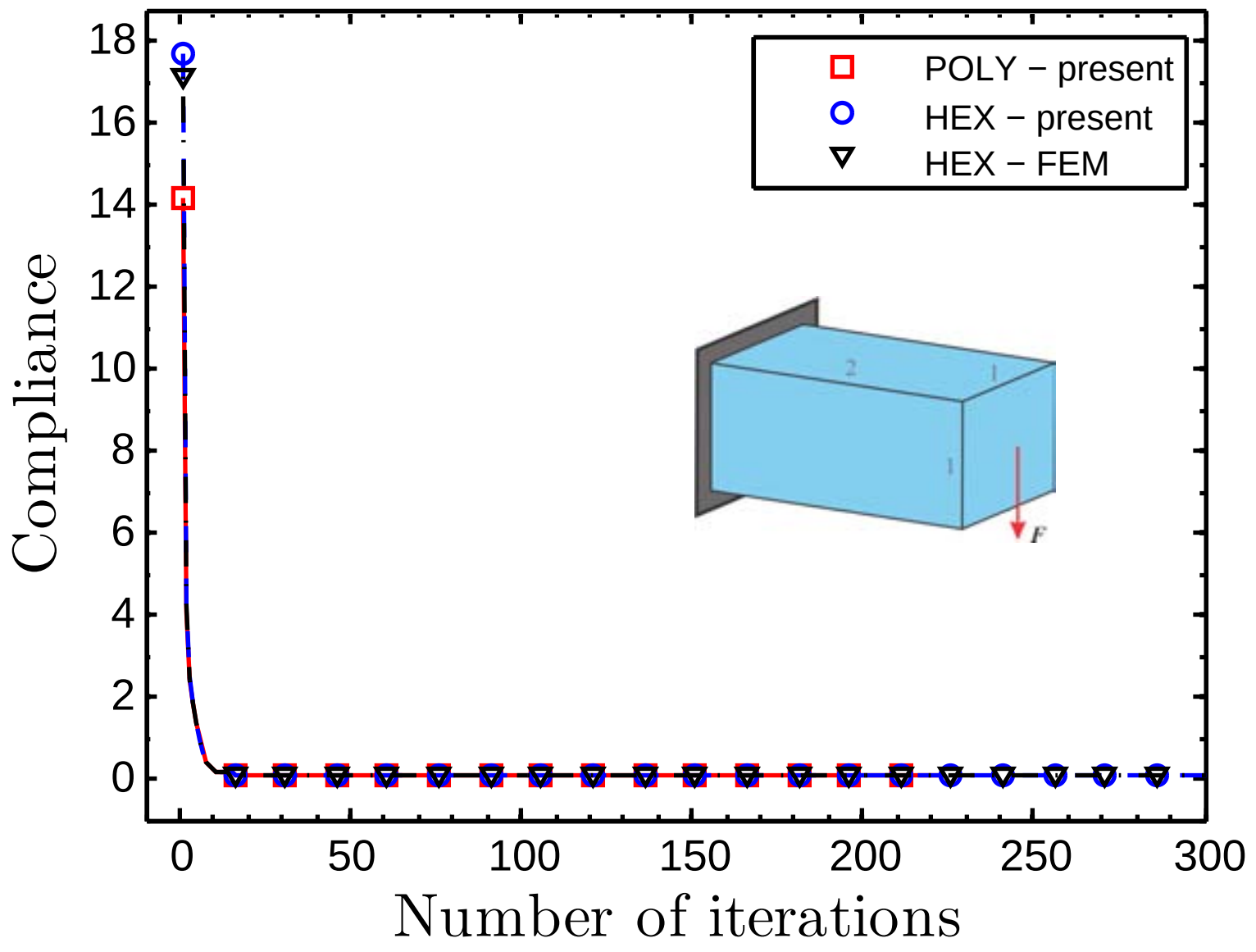

**Figure 7:** Convergence history for the cantilever beam problem.

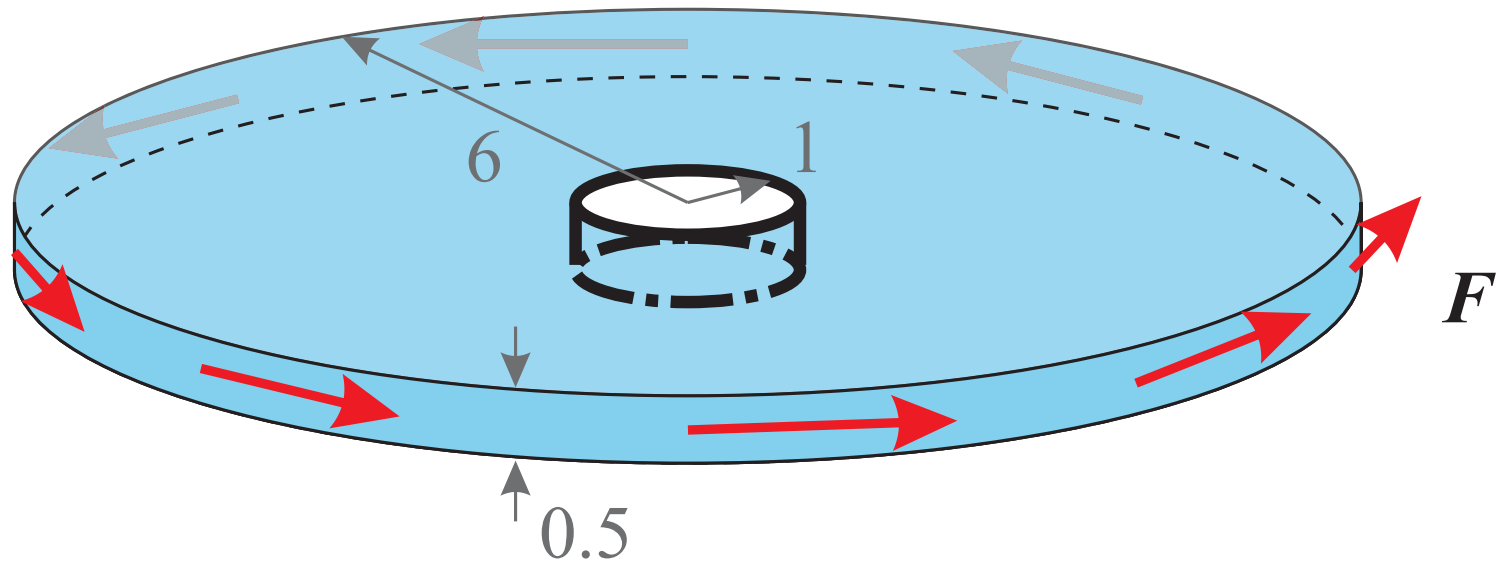

**Figure 8:** Problem description for thin disc. Eight equidistant shear loads are applied along the circumference of the disc and all the nodes along the central cylindrical hole are fixed.

According to Michell layout theory [30, 23], an optimum structural layout is one in which the tension and compression members meet orthogonally. Such a set of orthogonal curves are known as Hencky nets [23]. The tension and compression members in topology optimization solutions should adhere to this principle. The converged topology for the shear loaded thin disc problem, obtained from our algorithm, is shown in Fig. 9 (resembles a flower). The members of the structure intersect nearly at right angles, even for a coarse polyhedral mesh, indicating that mesh bias is alleviated with polyhedral elements. The compliance of the final topology is 0.5850.

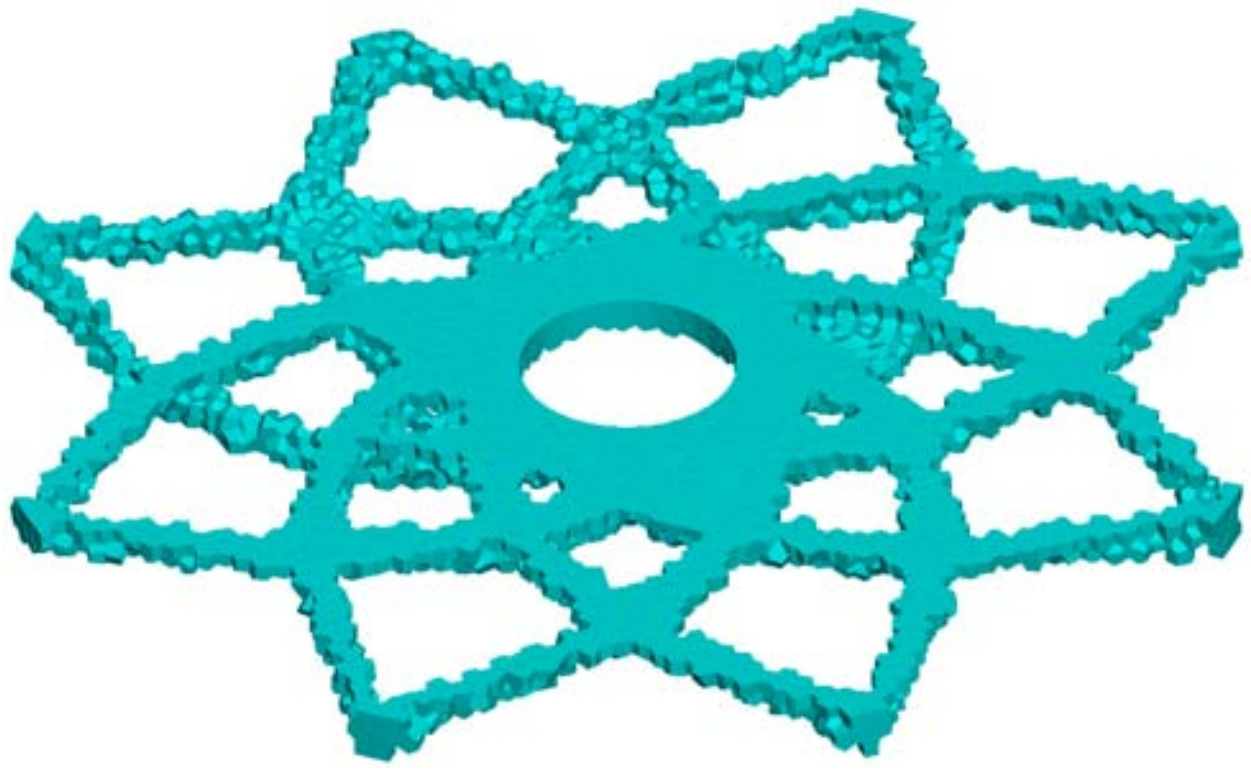

**Figure 9:** Final topology for shear loaded thin disc on a 10,000 element, 55,810 nodes polyhedral mesh. The average number of vertices per element is, $\mu = 20.53$, with standard deviation of, $\sigma = 3.71$, ($J = 0.5850$).

### *5.3. Hollow cylinder under torsional load*

We now look at the torsionally loaded hollow cylinder. The design domain is in the shape of a hollow cylinder of thickness 0.1 units, height 4 units and outer diameter of 1 unit (c.f. Fig. 10(a)). Four equidistant nodes along the bottom face are fixed and four tangential point loads are applied to corresponding nodes on the top face, effectively acting as a torsional load. A filter radius of 3% of the height of the cylinder and a volume fraction of 0.3 are prescribed. The problem is solved on three sets of meshes - two tetrahedral meshes of 9,977 elements (3,349 nodes); 451,584 elements (85,320 nodes) and a polyhedral mesh of 10,000 elements (79,925 nodes).

Using the polyhedral mesh, our optimization algorithm yields an elegant spiraling structure with nearly orthogonally oriented crossing members (Fig. 10(c)). For a tetrahedral mesh, with a similar number of elements as the polyhedral mesh (9,977), we still obtain a similar spiraling structure (c.f. Fig. 10(b)), but the members orientation is affected by the mesh geometry and the intersecting members are not fully orthogonal. This might be because although the tetrahedral mesh has similar number of elements as polyhedral mesh, the degrees of freedom in tetrahedral mesh is far less than the polyhedral mesh. So we solved the hollow cylinder problem on a fine tetrahedral mesh, such that the number of nodes in the mesh are comparable to that of the polyhedral mesh. Note that, in the fine tetrahedral mesh, the number of elements has risen to 451,584. Although the fine tetrahedral mesh (Fig. 10(d)) rectifies the lack of member orthogonality in the optimization result, the fine mesh considerably increases the computational cost associated with operations, such as creation and storage of the filter matrix, compared to the polyhedral and coarse tetrahedral meshes. If other length scale control approaches are adopted, for example enforcing perimeter constraint, in place of filters, the the computational cost associated with fine tetrahedral meshes can be considerably reduced. But the downside to this approach on tetrahedral meshes is that single node connections may arise in the designs. In terms of the cost associated with solving the governing elasticity problem, polyhedral meshes are marginally expensive than tetrahedral meshes of comparable total degrees of freedom. This is due to that fact that polyhedrons on an average have higher number of vertices (approximately 23 for our meshes) than tetrahedrons. Finally, note that the lower compliance for the optimization result on a coarse tetrahedral mesh (Fig. 10(b)) can be attributed to the fact that tetrahedral meshes experience artificial stiffness due to shear locking phenomenon, which reduces with mesh refinement.

We would like to point out that the element stiffness matrix obtained using the linear tetrahedral finite element approach is identical to the one obtained using the current approach. The reason being that, in the present method, for elements in the shape of a tetrahedron, the contribution of the stability term is zero, because the finite dimensional space of smooth functions on element, $\boldsymbol{\mathcal{W}}$, is identical to the space of linear deformations, $\boldsymbol{\mathcal{P}}$. So, the only contribution to the element stiffness matrix comes from the consistency term which is the same as the one obtained from finite element analysis. Thus, the final topologies obtained from the current approach and the linear tetrahedral finite element approach should be identical, along with the convergence history. Our results are in agreement with the above statement.

### *5.4. Hook domain under line load*

For the final compliance minimization problem, we investigate the hook domain subjected to a uniformly distributed line load along the negative $\boldsymbol{x}_{(3)}$-direction (Fig. 11(a)). A volume fraction of 0.1 is prescribed and a linear filter with radius equal to 2% of the maximum domain dimension is used. Using symmetry, we optimize only half the hook domain. The polyhedral mesh contains 10,000 elements (67,893 nodes). The converged topology, Fig. 11(b), has a compliance of 7.0484 and resembles the structure of a fan. The two dimensional version of the problem [44] has similar member orientations as our current three-dimensional result.

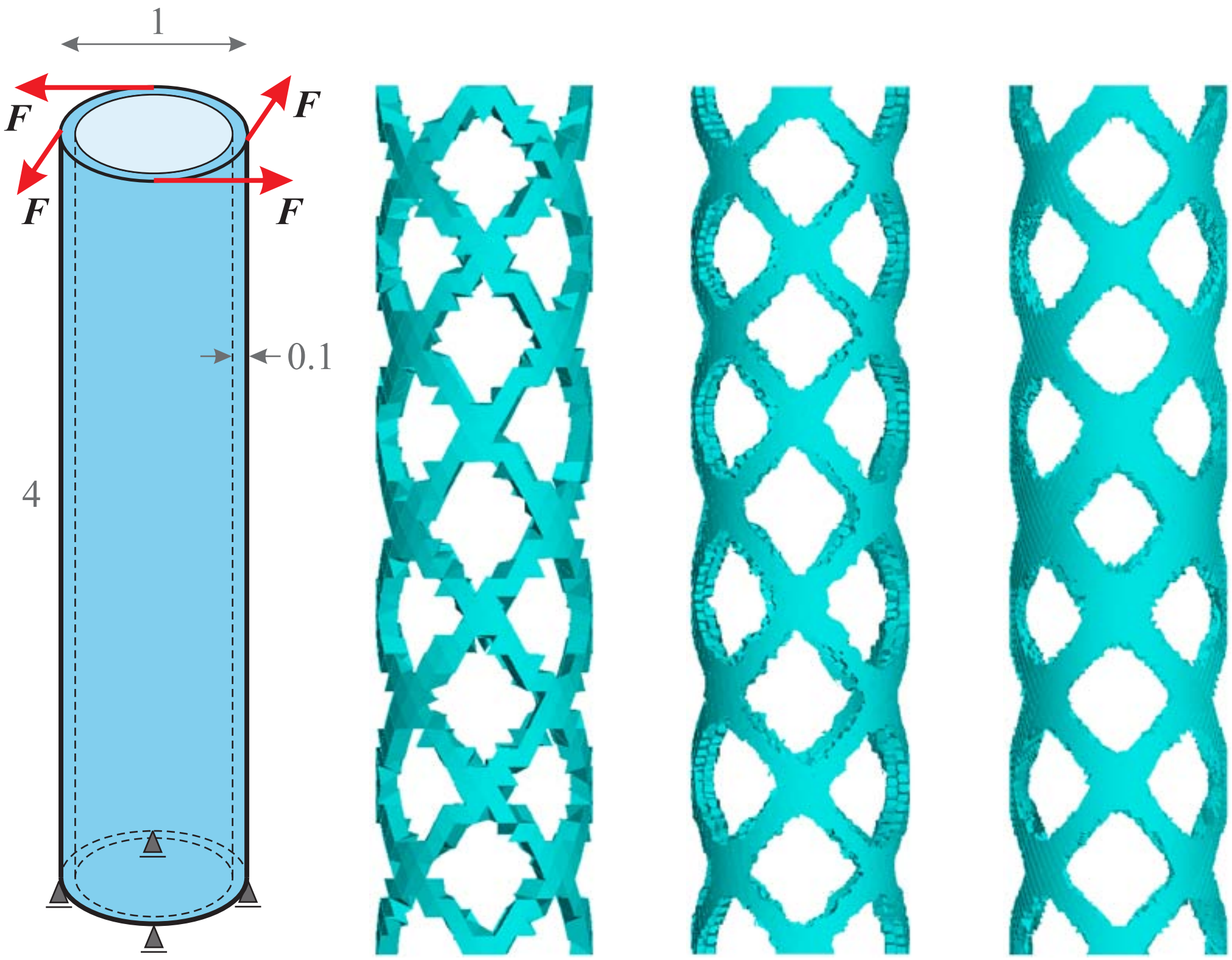

**Figure 10:** Hollow cylinder under torsional load. (a) Problem description. Converged topologies for (b) Linear tetrahedral mesh of 9,977 elements, 3,349 nodes ($J = 1.1397$); (c) Polyhedral mesh of 10,000 elements, 79,925 nodes. On average, polyhedral elements have, $\mu = 22.57$, vertices with standard derivation of, $\sigma = 2.88$, ($J = 1.6005$); (d) Linear tetrahedral mesh of 451,584 elements, 85,320 nodes ($J = 1.2064$).

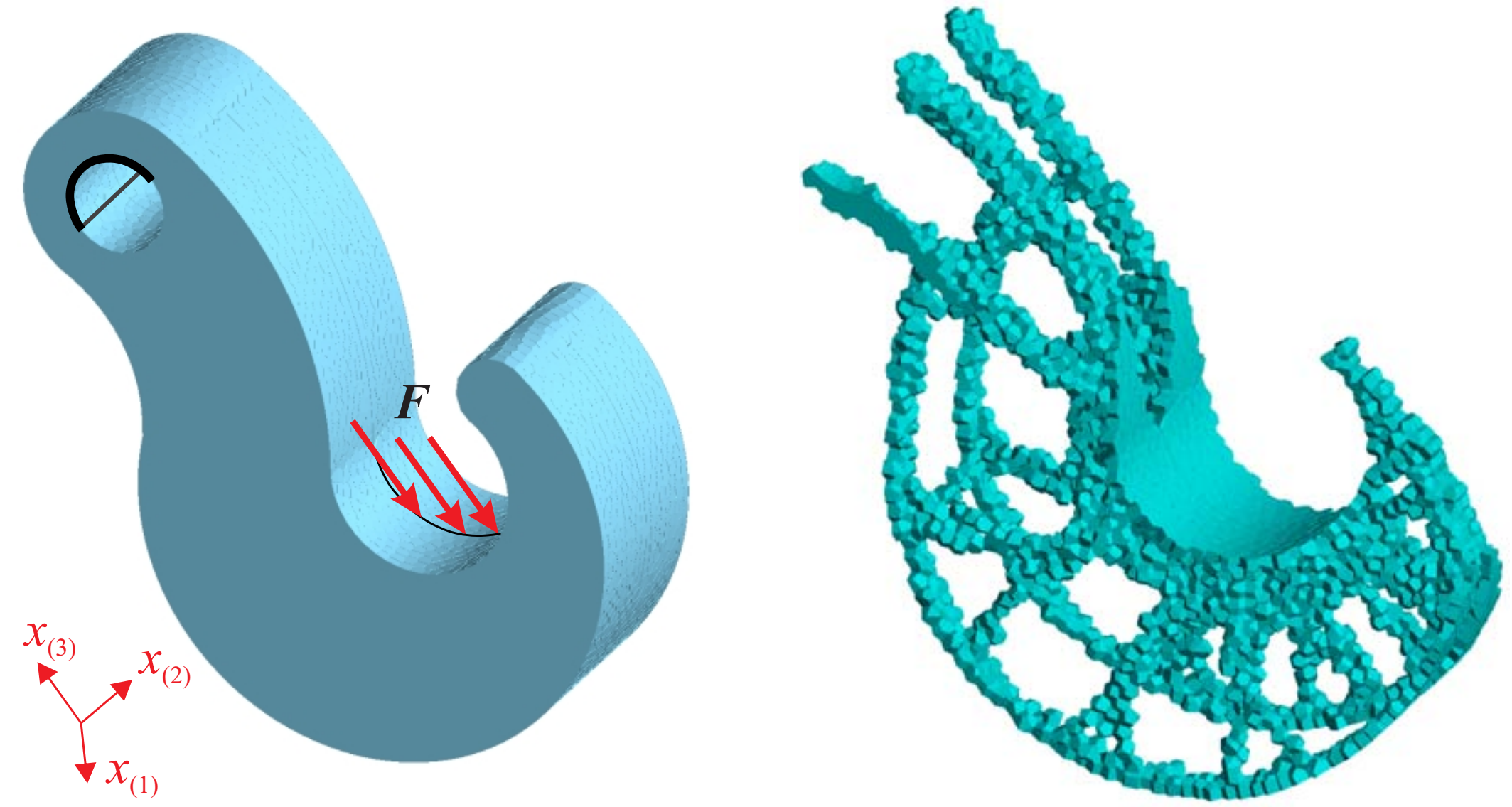

**Figure 11:** Hook under line load. (a) Problem description. All the nodes along the top half of the upper cylindrical hole are fixed and line load is applied along the circular arc in the negative $\boldsymbol{x}_{(3)}$-direction, as indicated in the figure. The domain is discretized using 10,000 polyhedral elements containing 67,893 nodes. The average number of vertices per polyhedron is, $\mu = 23.97$, with standard deviation, $\sigma = 4.19$. (b) Converged topology ($J = 7.0484$).

Next, we study the effect of filters on the optimization results for polyhedral meshes. Keeping all the parameters the same as before, we solve the hook problem without using any filter. As expected, without any length-scale control, more structural members, including some thin ones, appear in the solution (Fig. 12(b)). It is interesting to note that even without filtering no single node connections were observed in the design.

### *5.5. Displacement inverter*

Apart from compliance minimization, we also investigate a compliant mechanism problem. First, we explore the displacement inverter problem (Fig. 13). The domain is of dimension $1 \times 1 \times 1$ and is fixed at the bottom four corners. The objective of optimization is to maximize the output displacement $u_{out}$ on a workpiece modeled by a spring of stiffness $k_{out}$. The input and output spring stiffnesses, $k_{in}$ and $k_{out}$, are taken to be the same as the components of the global stiffness matrix at the degrees of freedom corresponding to the input and output nodes. Taking advantage of symmetry, only a quarter of the domain is optimized and is discretized using 10,000 polyhedral elements (58,785 nodes). A volume fraction of 0.05 is prescribed. Our final topology (Fig. 14) is similar to the ones available in the literature [1]. Visually, the three-dimensional optimization result (Fig. 14) resembles to the two-dimensional optimization result [10].

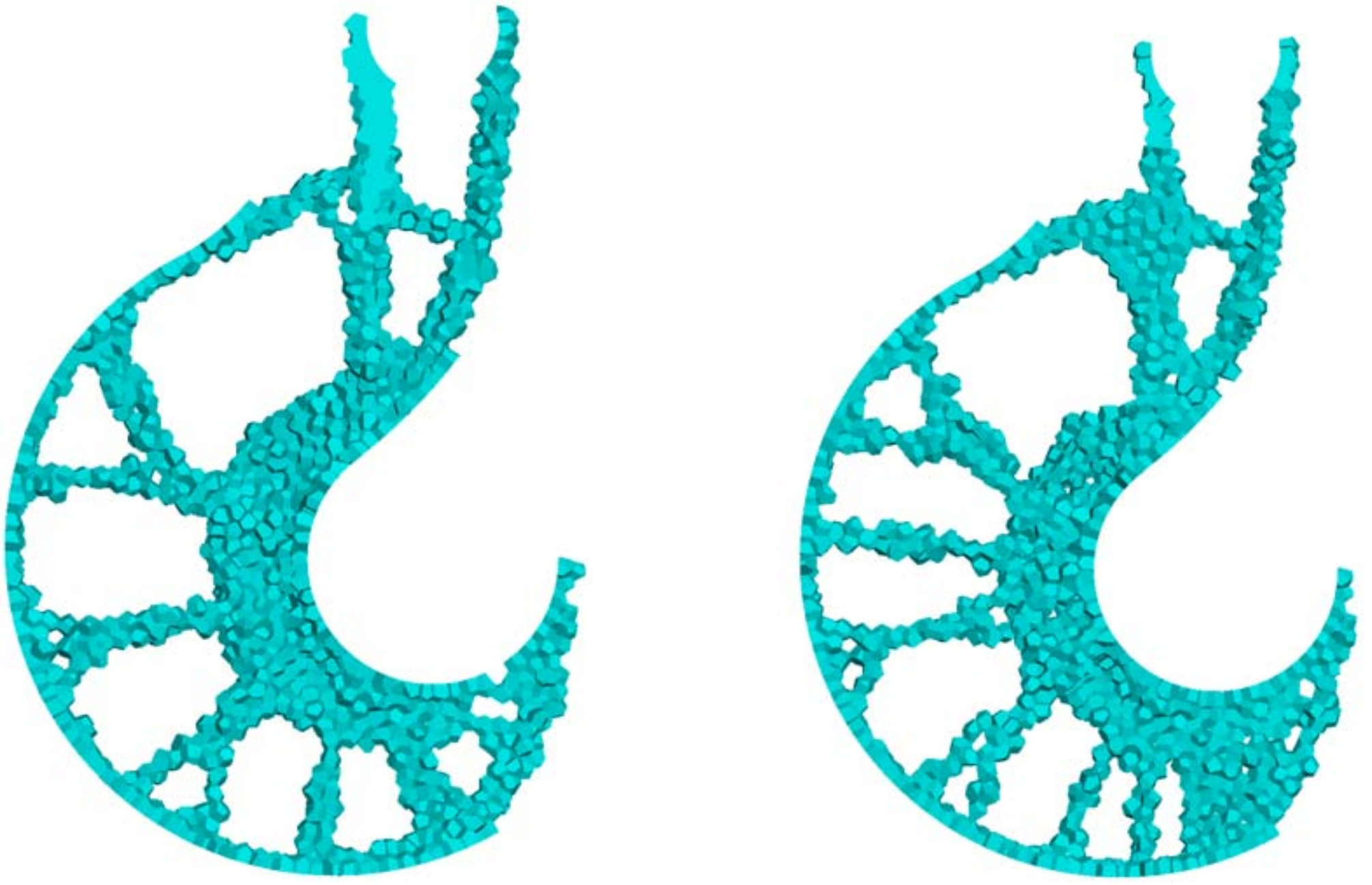

**Figure 12:** Study of the influence of filters on polyhedral topology optimization for hook problem of Fig. 11(a). (a) With filter (front view of Fig. 11(b)). (b) Without filter.

### *5.6. Gripper*

A gripper, as the name suggests, is a complaint mechanism suitable for gripping objects. The input actuator, modeled as a spring with stiffness $k_{in}$, and a force $f_{in}$, applies a horizontal load as indicated in Fig. 15(a) and the goal of optimization is to maximize the vertical output displacement $u_{out}$ on a workpiece modeled by a spring of stiffness $k_{out}$. The problem dimensions and boundary conditions are indicated in Fig. 15(a). Nodes on the top and bottom section on the right face, indicated by gray color, are fixed. Using the passive element concept the orange box is modeled as void. Before each optimization cycle, the elements lying inside the orange box are identified and are assigned the minimum density corresponding to voids. The spring stiffnesses, $k_{in}$ and $k_{out}$ are taken to be the same as the components of the global stiffness matrix at the degrees of freedom corresponding to the input and output nodes. Due to symmetry, only a quarter of the domain is optimized and is discretized using 10,000 polyhedral elements (58,785 nodes) with a prescribed volume fraction of 0.10. The converged topology (Fig. 15) for three-dimensional optimization is a clear extension of the two-dimensional results available in reference [10]. The combination of gripping jaws and hinge mimicking narrow necks, along the middle of the design, resemble a pair of scissors. The convergence history for the gripper problem, shown in Fig. 16, indicates stable convergence.

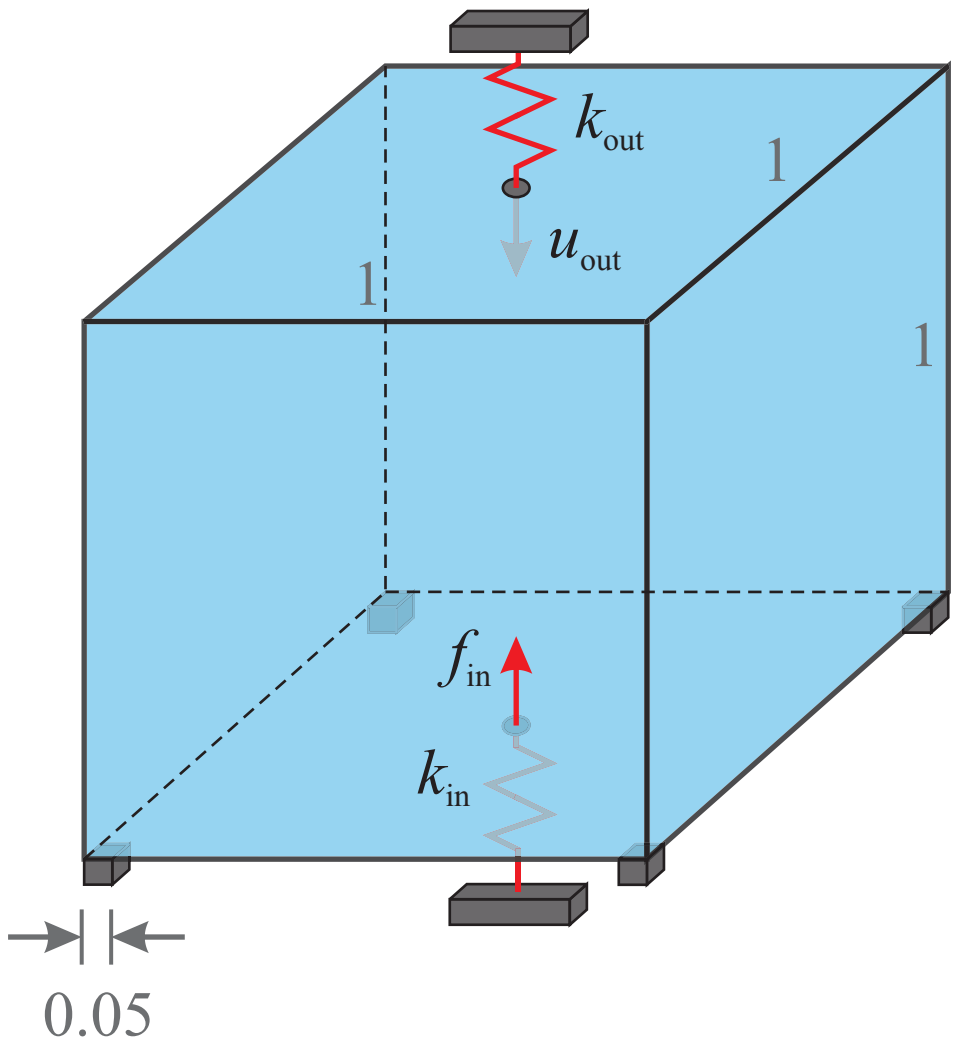

**Figure 13:** Displacement inverter.

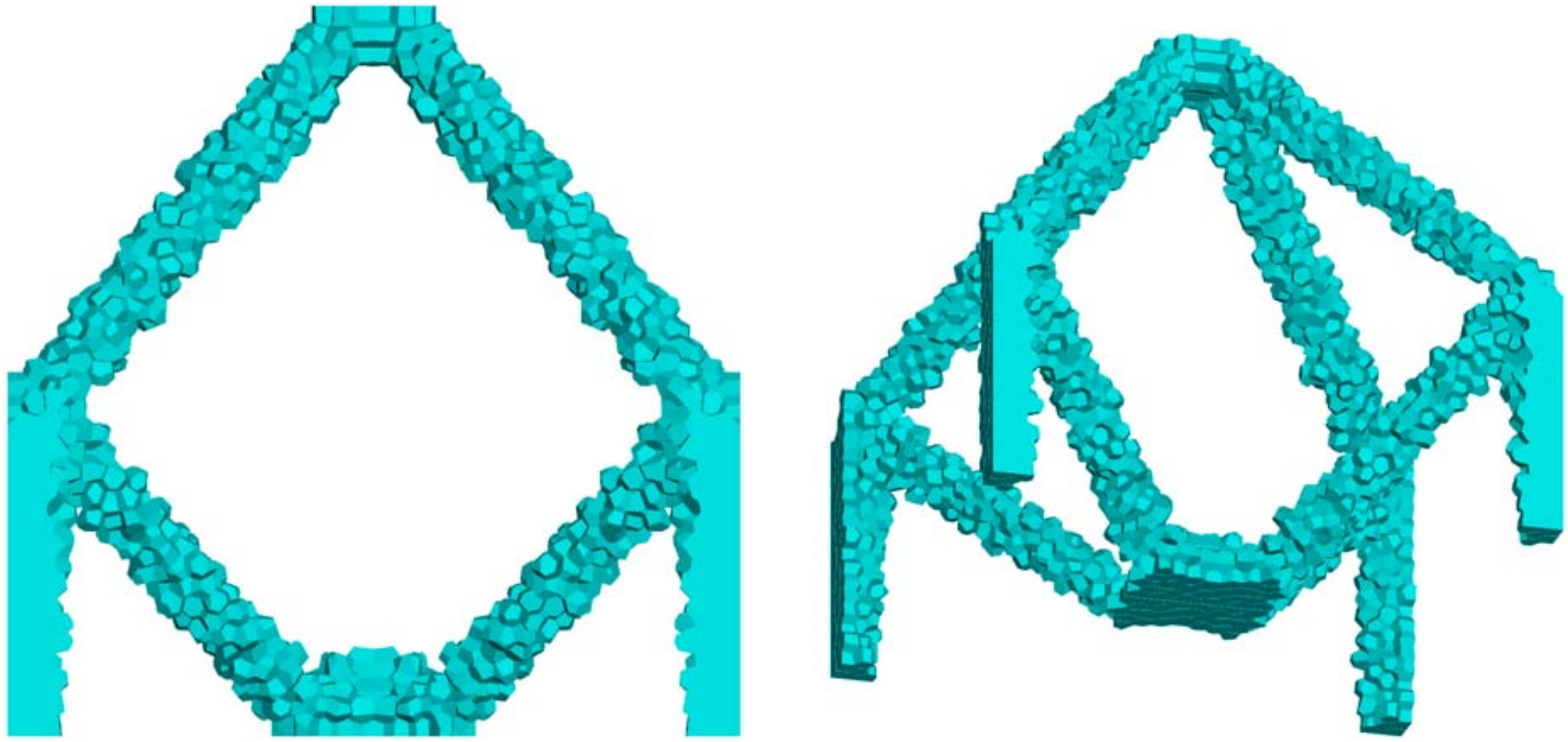

**Figure 14:** Converged topology for the displacement inverter problem. Due to symmetry, only a quarter of the domain is optimized and is discretized using 10,000 polyhedral elements containing 58,785 nodes. The average number of vertices per polyhedron is, $\mu = 22.98$, with standard deviation of, $\sigma = 3.75$.

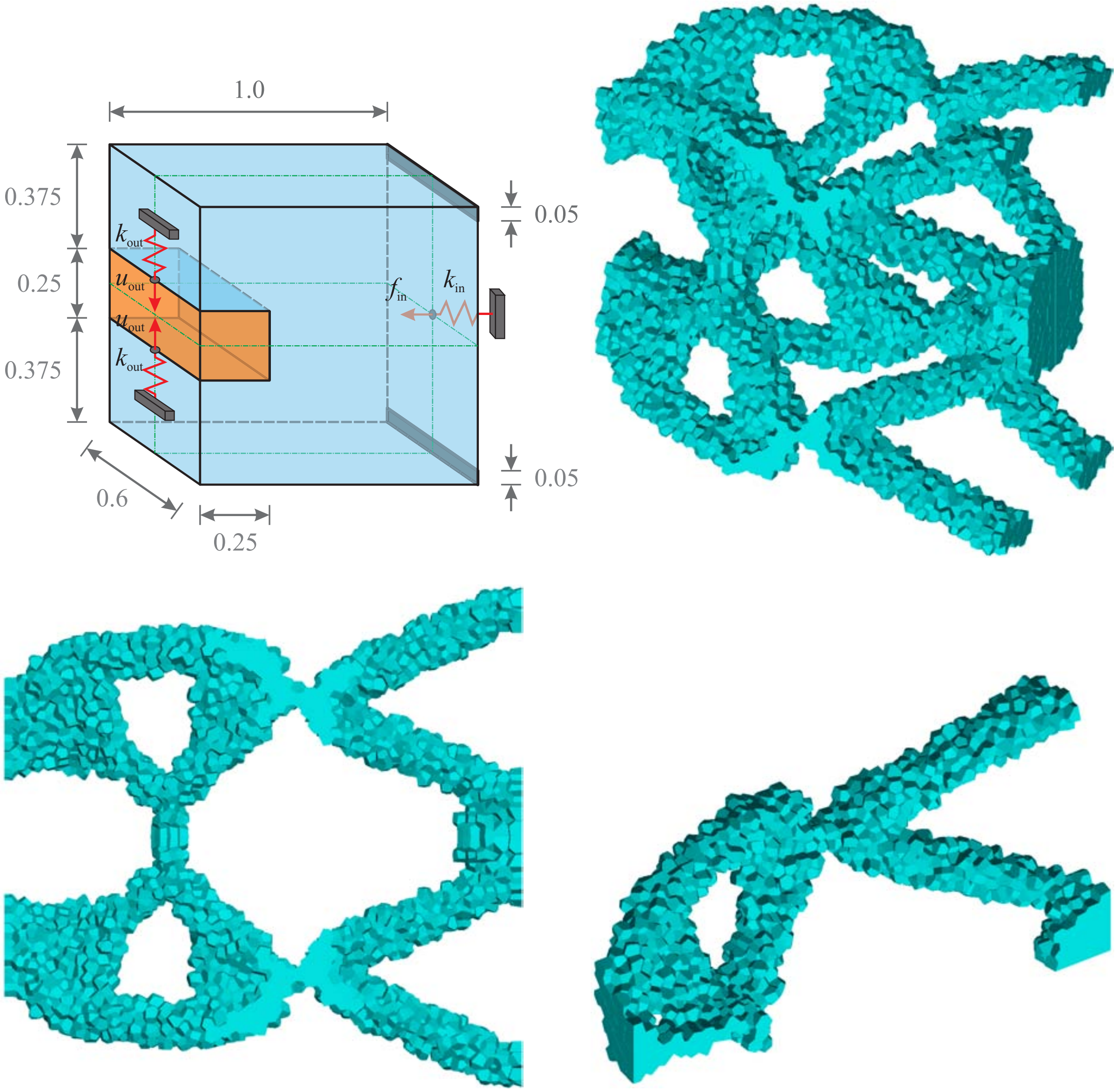

**Figure 15:** Topology optimization design of gripper. A quarter of the problem is solved with a polyhedral mesh of 10,000 elements and 59,194 nodes. The average number of vertices per polyhedron is, $\mu = 23.11$, with standard deviation, $\sigma = 3.85$. Orange region is modeled as voids using passive element concept. (a) Problem description. (b) Complete design. (c) Front view. (d) Quarter section view.

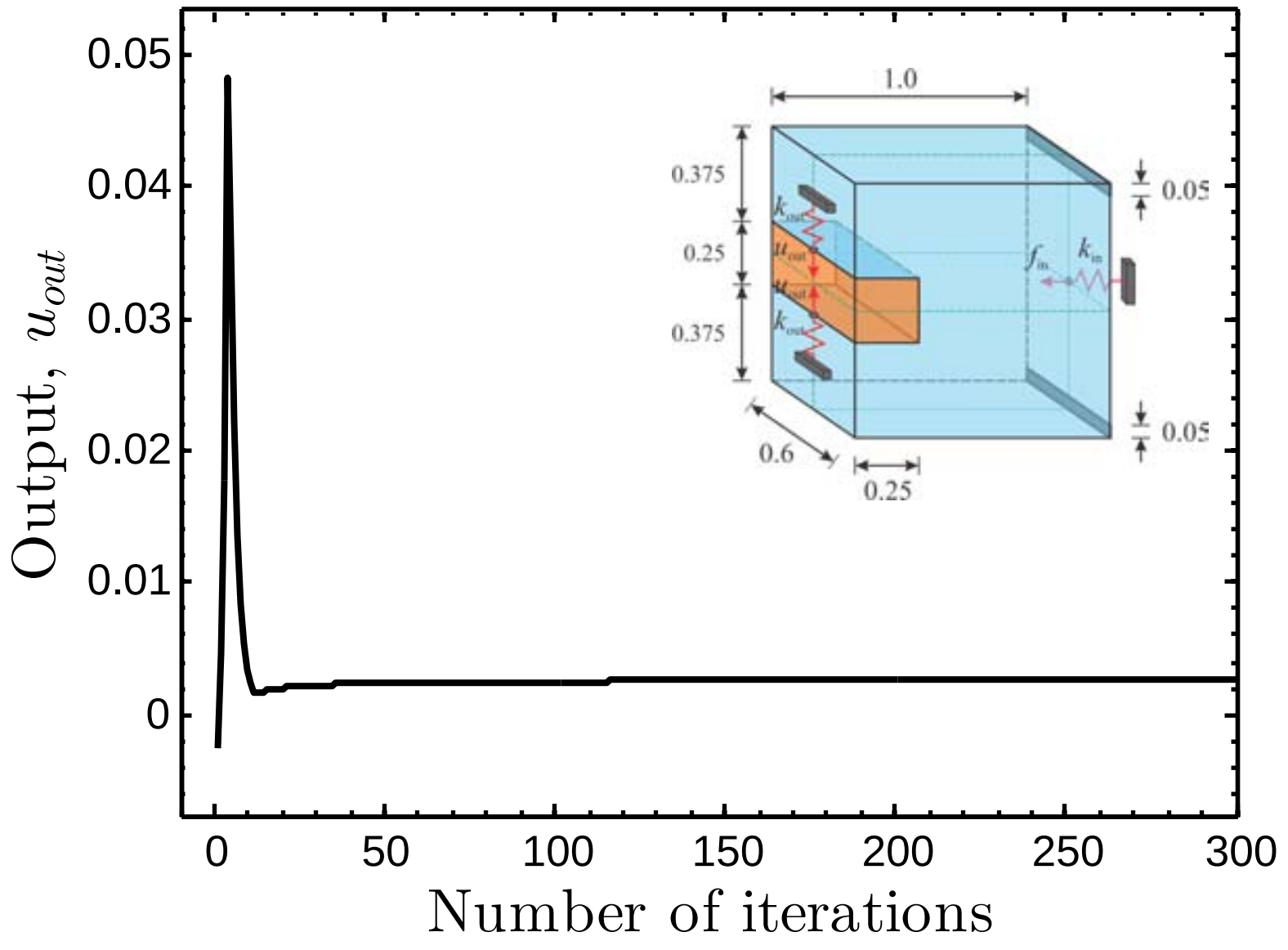

**Figure 16:** Convergence history for the gripper problem.

## 6. Concluding remarks

In this work, we explore polyhedral meshes for three-dimensional topology optimization. We use the Virtual Element Method to numerically solve the governing elasticity problem. VEM addresses some of the challenges facing finite element approaches, such as expensive shape function and their derivatives computation, accurate numerical evaluation of the weak form and sensitivity to degenerate elements. Unlike other Galerkin approaches, VEM does not require explicit computation of the shape functions. Rather the stiffness matrix is constructed directly using a projection map which is defined to extract the constant strain deformation modes, thereby guaranteeing the passing of patch test.

We investigate the topology optimization of compliance minimization and compliant mechanism problems using polyhedrons which are generated using centroidal Voronoi tessellation (CVT) based meshing approach. Compliance minimization problems such as the cantilever beam, shear loaded disc, hollow cylinder under torsion, hook problem subjected to line load; and compliant mechanism problems such as the displacement inverter and gripper are explored. As indicated by our thin disc and hollow cylinder problem solutions, relatively coarse arbitrary polyhedral elements capture member orthogonality and alleviate mesh bias designs. Numerical anomalies of single node connections and checkerboarding, prevalent in topology optimization, are naturally alleviated by polyhedrons due to their intrinsic geometry. Current work of three-dimensional topology optimization with polyhedrons paves the way to future applications in the field of multiphysics designs [35, 13] and biomedical engineering [40]. For large scale topology optimization problems, a multiresolution approach [20, 31] may be adopted, where in the elasticity and optimization problem are solved on a coarse and fine mesh, respectively.

## Nomenclature

| | |
|---|---|
| $a$, $a_h$ | continuous and discrete bilinear form |
| $f$ | continuous load linear form |
| $\boldsymbol{t}$ | surface tractions |
| $\boldsymbol{C}$ | elasticity tensor |
| $\Omega$ | topology optimization working domain |
| $\boldsymbol{\mathcal{V}}$, $\boldsymbol{\mathcal{V}}_h$ | continuous and conforming discrete solution space |
| $J$ | objective function |
| $\rho$, $\rho_e$ | continuous density function and discrete element density |
| $V_f$ | prescribed volume fraction |
| $\lvert\cdot\rvert$ | measure (area or volume) of a set. Also, Euclidean norm of vector |
| $\boldsymbol{\mathcal{W}}$ | space of finite dimensional smooth functions over element $e$ |
| $a^e$, $a_h^e$ | restriction of $a$, $a_h$ to $\boldsymbol{\mathcal{W}}$ (element contribution) |
| $s^e$ | approximation of bilinear form corresponding to higher-order deformation modes |
| $\varphi_i$ | generic barycentric coordinates |
| $\boldsymbol{\varphi}_i$ | canonical basis function |
| $\overline{\boldsymbol{v}}$ | mean of values of $\boldsymbol{v}$ sampled over the vertices of element $e$ |
| $\langle\boldsymbol{v}\rangle$ | volume average of $\boldsymbol{v}$ |
| $\boldsymbol{\mathcal{P}}$ | space of linear deformation mode |
| $\boldsymbol{p}_i$ | bases spanning space $\boldsymbol{\mathcal{P}}$ |
| $\pi_{\boldsymbol{\mathcal{P}}}$ | projection map to extract $\boldsymbol{\mathcal{P}}$, $\pi_{\boldsymbol{\mathcal{P}}} : \boldsymbol{\mathcal{W}} \to \boldsymbol{\mathcal{P}}$ |
| $\boldsymbol{P}_{\boldsymbol{\mathcal{P}}}$ | matrix representation of $\pi_{\boldsymbol{\mathcal{P}}}$ |
| $\boldsymbol{W}_{\boldsymbol{\mathcal{P}}}$ | matrices containing surface integration quantities |
| $\boldsymbol{N}_{\boldsymbol{\mathcal{P}}}$ | matrices containing rearranged nodal coordinates of vertices |
| $\boldsymbol{q}_i$ | vector of facial surface integration of barycentric coordinates |
| $\boldsymbol{K}_e$ | element stiffness matrix |
| $\boldsymbol{D}$ | material matrix which is function of elasticity tensor $\boldsymbol{C}$ |
| $\boldsymbol{I}$ | identity matrix |
| $\alpha^e$ | positive scaling coefficient for stability term $s^e$ |
| $\boldsymbol{\mathcal{S}}$ | set of distinct seeds to generate Voronoi cells |
| $E$ | Young's modulus |
| $\nu$ | Poisson's ratio |
| $\boldsymbol{K}$ | global stiffness matrix |
| $\boldsymbol{F}$ | global force vector |
| $\boldsymbol{U}$ | global displacement vector |
| $r_{min}$ | filter radius |

## Appendix A: Sensitivity analysis

The discretized version of (6) can be expressed as:

$$\begin{aligned} &\inf_{\rho} \; J = \boldsymbol{P}^T \boldsymbol{U} \\ &\text{subject to:} \quad \boldsymbol{K}(\rho)\boldsymbol{U} = \boldsymbol{F}, \qquad \textstyle\sum_e |e|\rho_e \leq V_f|\Omega| \end{aligned} \tag{39}$$

where $\boldsymbol{K}$, $\boldsymbol{U}$ and $\boldsymbol{F}$ are the global stiffness matrix, global nodal displacement vector and global nodal force vector, respectively. The vector $\boldsymbol{P}$ represents the global force vector $\boldsymbol{F}$ for compliance minimization problem because the objective is to minimize external work. For displacement inverter and gripper problem, $\boldsymbol{P}$ is a vector with all zeros except at locations corresponding to the output degree of freedom, where it is unity.

We use gradient-based optimization algorithm (OC) for solving the discrete problem (39), which requires computation of the gradient of the objective function $J$. Using the adjoint method [10], the sensitivity of $J$ with respect to the design variable (element density, $\rho_e$), is given by:

$$\frac{\partial J}{\partial \rho_e} = -p(1-\epsilon)\rho_e^{p-1}\boldsymbol{\lambda}_e^T \boldsymbol{K}_e \boldsymbol{u}_e \tag{40}$$

For the compliance minimization problem (3), the vector $\boldsymbol{\lambda}_e$ is the same as the element displacement vector $\boldsymbol{u}_e$. For the displacement inverter problem (5), $\boldsymbol{\lambda}_e$ is the elemental component of $\boldsymbol{\lambda}$ that solves the adjoint system $\boldsymbol{K}\boldsymbol{\lambda} = \boldsymbol{P}$. Also, the sensitivity of the volume constraint (shown in (39)) with respect to element density is $|e|$.

## Appendix B: Projection method - Filtering

In topology optimization, filters are used to enforce a length-scale in the problem and to ensure mesh-independency. We use a linear filter which assigns a weighted average of the nearby elemental densities to each element. Thus, the projected element density of an element, $\rho_e$, is written as:

$$\rho_e = \frac{\sum_{j \in N_e} w_{ej}\rho_j}{\sum_{j \in N_e} w_{ej}} \tag{41}$$

Here, $N_e$ is the set of all elements whose centers lie within a distance of $r_{min}$ from the center of the element $e$ under consideration. The linear weights $w_{ej}$ are:

$$w_{ej} = \max\left(0, \frac{r_{min} - r_{ej}}{r_{min}}\right) \tag{42}$$

where $r_{min}$ and $r_{ej}$ are the enforced minimum member size and the distance between centroids of elements $e$ and $j$, respectively.

To compute the sensitivity of the objective function $J$ with respect to the independent design variables $\rho_j$, we use the chain rule:

$$\frac{\partial J}{\partial \rho_j} = \frac{\partial J}{\partial \rho_e}\frac{\partial \rho_e}{\partial \rho_j}, \qquad \text{where } \frac{\partial \rho_e}{\partial \rho_j} = \frac{w_{ej}}{\sum_{j \in N_e} w_{ej}} \tag{43}$$

Other filters, such as sensitivity filters may also be used.

## References

[1] O. Amir and O. Sigmund. On reducing computational effort in topology optimization: how far can we go? *Structural and Multidisciplinary Optimization*, 44:25–29, 2011.

[2] M. Arroyo and M. Ortiz. Local maximum-entropy approximation schemes: a seamless bridge between finite elements and meshfree methods. *International Journal for Numerical Methods in Engineering*, 65(13):2167–2202, 2006.

[3] F. Aurenhammer. Voronoi diagrams - a survey of a fundamental geometric data structure. *ACM Computing Surveys*, 23(3):345–405, 1991.

[4] L. Beirão Da Veiga. A mimetic discretization method for linear elasticity. *ESAIM: Mathematical Modelling and Numerical Analysis*, 44(2):231–250, 2010.

[5] L. Beirão Da Veiga, F. Brezzi, A. Cangiani, G. Manzini, L. D. Marini, and A. Russo. Basic principles of virtual element methods. *Mathematical Models and Methods in Applied Sciences*, 23(1):199–214, 2013.

[6] L. Beirão Da Veiga, F. Brezzi, and L. D. Marini. Virtual Elements for linear elasticity problems. *SIAM Journal on Numerical Analysis*, 2013.

[7] L. Beirão Da Veiga, K. Lipnikov, and G. Manzini. Error analysis for a mimetic discretization of the steady stokes problem on polyhedral meshes. *SIAM Journal on Numerical Analysis*, 48(4):1419–1443, 2010.

[8] V. V. Belikov, V. D. Ivanov, V. K. Kontorovich, S. A. Korytnik, and A. Yu. Semenov. The non-Sibsonian interpolation: A new method of interpolation of the values of a function on an arbitrary set of points. *Computational Mathematics and Mathematical Physics*, 37(1):9–15, 1997.

[9] M. P. Bendsøe. Optimal shape design as a material distribution problem. *Structural Optimization*, 1:193–202, 1989.

[10] M. P. Bendsøe and O. Sigmund. *Topology Optimization - Theory, Methods and Applications.* Springer: New York, 2003.

[11] J. Bishop. A displacement-based finite element formulation for general polyhedra using harmonic shape functions. *International Journal for Numerical Methods in Engineering*, 2013. Accepted for publication.

[12] F. Brezzi, K. Lipnikov, M. Shashkov, and V. Simoncini. A new discretization methodology for diffusion problems on generalized polyhedral meshes. *Computer Methods in Applied Mechanics and Engineering*, 196(37-40 SPEC. ISS.):3682–3692, 2007.

[13] R. C. Carbonari, E. C. N. Silva, and Paulino. G. H. Topology optimization design of functionally graded bimorph-type piezoelectric actuators. *Smart Materials & Structures*, 16(6):2605–2620, 2007.

[14] N. H. Christ, R. Friedberg, and Lee T. D. Weights of links and plaquettes in a random lattice. *Nuclear Physics B*, 210(3):337–346, 1982.

[15] Q. Du and D. Wang. The optimal centroidal Voronoi tessellations and the gersho's conjecture in the three-dimensional space. *Computers and Mathematics with Applications*, 49(9-10):1355–1373, 2005.

[16] P. Dvorak. New element lops time off CFD simulations. *Machine Desig*, 78(5):154–155, 2006.

[17] C. Fleury and V. Braibant. Structural optimization: a new dual mehod using mixed variables. *International Journal for Numerical Methods in Engineering*, 23(3):409–428, 1986.

[18] M. S. Floater. Mean value coordinates. *Computer Aided Geometric Design*, 20(1):19–27, 2003.

[19] M. S. Floater, G. Kòs, and M. Reimers. Mean value coordinates in 3D. *Computer Aided Geometric Design*, 22(7):623–631, 2005.

[20] A. L. Gain. *Polytope-based topology optimization using a mimetic-inspired method.* PhD thesis, University of Illinois at Urbana-Champaign, 2013.

[21] A. L. Gain and G. H. Paulino. Phase-field based topology optimization with polygonal elements: A finite volume approach for the evolution equation. *Structural and Multidisciplinary Optimization*, 46(3):327–342, 2012.

[22] A. L. Gain, C. Talischi, and G. H. Paulino. On the virtual element method for three-dimensional elasticity problems on arbitrary polyhedral meshes. 2013. *Submitted.*

[23] W. S. Hemp. *Optimum Structures.* Clarendon Press, Oxford, 1973.

[24] H. Hiyoshi and K. Sugihara. Two generalizations of an interpolant based on Voronoi diagrams. *International Journal of Shape Modeling*, 5(2):219–231, 1999.

[25] K. Hormann and N. Sukumar. Maximum entropy coordinates for arbitrary polytopes. In *Eurographics Symposium on Geometry Processing*, volume 27, pages 1513–1520, 2008.

[26] P. Joshi, M. Meyer, T. Derose, B. Green, and T. Sanocki. Harmonic coordinates for character articulation. *ACM Transactions on Graphics*, 26(3), 2007. Art. no. 1276466.

[27] S. Lloyd. Least squares quantization in PCM. *IEEE Transactions on Information Theory*, 28(2):129–137, 1982.

[28] E. A. Malsch, J. J. Lin, and G. Dasgupta. Smooth two dimensional interpolants: a recipe for all polygons. *Journal of Graphics Tools*, 10:2, 2005.

[29] S. Martin, P. Kaufmann, M. Botsch, M. Wicke, and M. Gross. Polyhedral finite elements using harmonic basis functions. *Computer Graphics Forum*, 27(5):1521–1529, 2008.

[30] A. G. M. Michell. The limits of economy of material in frame-structures. *Philosophical Magazine*, 8(47):589–597, 1904.

[31] T. H. Nguyen, G. H. Paulino, J. Song, and C. H. Le. A computational paradigm for multiresolution topology optimization (MTOP). *Structural and Multidisciplinary Optimization*, 41(4):525–539, 2010.

[32] P.-O. Persson and G. Strang. A simple mesh generator in MATLAB. *SIAM Review*, 46(2):329–345, 2004.

[33] G. I. N. Rozvany, M. Zhou, and T. Birker. Generalized shape optimization without homogenization. *Structural and Multidisciplinary Optimization*, 4(3-4):250–252, 1992.

[34] R. Sibson. A vector identity for the Dirichlet tessellation. *Mathematical Proceedings of the Cambridge Philosophical Society*, 87:151–155, 1980.

[35] O. Sigmund. Design of multiphysics actuators using topology optimization - Part II: Two-material structures. *Computer Methods in Applied Mechanics and Engineering*, 190(49-50):6605–6627, 2001.

[36] N. Sukumar. Construction of polygonal interpolants: a maximum entropy approach. *International Journal of Numerical Methods in Engineering*, 61(12):2159–2181, 2004.

[37] N. Sukumar and E. A. Malsch. Recent advances in the construction of polygonal finite element interpolations. *Archives of Computational Methods in Engineering*, 13(1):129–163, 2006.

[38] N. Sukumar, B. Moran, A. Y. Semenov, and V. V. Belikov. Natural neighbor Galerkin methods. *International Journal for Numerical Methods in Engineering*, 50:1–27, 2001.

[39] N. Sukumar and A. Tabarraei. Conforming polygonal finite elements. *International Journal of Numerical Methods in Engineering*, 61(12):2045–2066, 2004.

[40] A. Sutradhar, G. H. Paulino, M. J. Miller, and T. H. Nguyen. Topology optimization for designing patient-specific large craniofacial segmental bone replacements. *Proceedings of the National Academy of Sciences*, 107(30):13222–13227, 2010.

[41] K. Svanberg. The method of moving asymptotoes - a new method for structural optimization. *International Journal for Numerical Methods in Engineering*, 24(2):359–373, 1987.

[42] C. Talischi, G. H. Paulino, A. Pereira, and I. F. M. Menezes. Polygonal finite elements for topology optimization: A unifying paradigm. *International Journal for Numerical Methods in Engineering*, 82:671–698, 2010.

[43] C. Talischi, G. H. Paulino, A. Pereira, and I. F. M. Menezes. PolyMesher: A general-purpose mesh generator for polygonal elements written in Matlab. *Structural and Multidisciplinary Optimization*, 45(3):309–328, 2012.

[44] C. Talischi, G. H. Paulino, A. Pereira, and I. F. M. Menezes. PolyTop: a Matlab implementation of a general topology optimization framework using unstructured polygonal finite element meshes. *Journal of Structural and Multidisciplinary Optimization*, 45(3):329–357, 2012.

[45] E. L. Wachspress. *A rational finite element basis.* Academic Press, New York, 1975.

[46] E. L. Wachspress. Rational bases for convex polyhedra. *Computers and Mathematics with Applications*, 59(6):1953–1956, 2010.

[47] J. Warren. Barycentric coordinates for convex polytopes. *Advances in Computational Mathematics*, 6(1):97–108, 1996.

[48] M. Yip, J. Mohle, and J. E. Bolander. Automated modeling of three-dimensional structural components using irregular lattices. *Computer Aided Civil and Infrastructure Engineering*, 20(6):393–407, 2005.